\documentclass[a4paper,12pt]{article} 
\usepackage{amsmath} 
\usepackage{amsopn} 
\usepackage{amssymb} 
\usepackage[mathscr]{eucal} 
\usepackage{theorem} 
\usepackage{enumerate} 

\theorembodyfont{\itshape} 
\theoremstyle{plain}
  \newtheorem{thm}{Theorem}[section] 
  \newtheorem{pro}[thm]{Proposition} 
   
  \newtheorem{cor}[thm]{Corollary} 

\theorembodyfont{\rmfamily} 
\theoremstyle{plain}
  \newtheorem{rem}{Remark}[section]

\renewcommand{\theequation}%
           {\thesection.\arabic{equation}}

\begin{document} 

\begin{center} 
{\LARGE The equations of Gauss, Codazzi and Ricci} 

\vspace{1mm} 

{\LARGE of surfaces in 4-dimensional space forms} 

\vspace{6mm} 

{\Large Naoya {\sc Ando}} 
\end{center} 

\vspace{3mm} 

\begin{quote} 
{\footnotesize \it Abstract} \ 
{\footnotesize 
Let $N$ be a Riemannian, neutral or Lorentzian $4$-dimensional space form. 
In this paper, 
the expressions of the equations of Gauss, Codazzi and Ricci 
of a space-like or time-like surface in $N$ 
given in \cite{AH} are naturally understood 
in terms of the induced connection (of the complexification) 
of the two-fold exterior power of the pull-back bundle on the surface. 
Moreover, based on such expressions, 
we characterize several classes of surfaces 
related to the covariant derivatives of the twistor lifts and so on.} 
\end{quote} 

\vspace{3mm} 

\section{Introduction}\label{sect:intro} 

\setcounter{equation}{0} 

Let $N$ be a Riemannian, neutral or Lorentzian $4$-dimensional space form. 
In this paper, we understand the equations of Gauss, Codazzi and Ricci 
of a space-like or time-like surface in $N$, 
in terms of the induced connection (of the complexification) 
of the two-fold exterior power of the pull-back bundle on the surface. 

Suppose that $N$ is Riemannian. 
Then the two-fold exterior power of the pull-back bundle 
on a surface in $N$ is 
decomposed into two orientable subbundles of rank $3$. 
Computing the curvature tensors of these subbundles, 
we will obtain the expressions of the equations of Gauss, Codazzi and Ricci 
of the surface given in \cite{AH} (Proposition~\ref{pro:GCRR}). 
The twistor spaces associated with the pull-back bundle are 
the unit sphere bundles of these subbundles. 
The twistor lifts of the surface are sections of these twistor spaces. 
Refer to \cite{bryant}, \cite{ES}, \cite{friedrich} for twistor spaces 
and twistor lifts. 
We will define the degeneracy of the twistor lifts 
and observe that the degeneracy of the twistor lifts just means that 
the surface satisfies the following: 
\begin{itemize} 
\item[{\rm (i)}]{the curvature $K$ of the induced metric of the surface is 
identically equal to 
the constant sectional curvature $L_0$ of $N$;} 
\item[{\rm (ii)}]{the normal connection $\nabla^{\perp}$ of the surface is 
flat} 
\end{itemize} 
(see Proposition~\ref{pro:WX+YZ=0}). 
Refer to \cite{AH}, \cite{chen}, \cite{DT}, \cite{enomoto}, \cite{hoffman}, 
\cite{MS} and \cite{thas} for surfaces with flat normal connection. 
We will obtain characterizations of a surface 
such that the twistor lifts are nondegenerate, 
by relations among functions in the above expressions related to 
the second fundamental form 
(Propositions~\ref{pro:nondegL=0}, \ref{pro:nondegLnot=0}). 
An isotropic minimal surface is characterized by 
horizontality of one of the twistor lifts (\cite{friedrich}). 
We will define the complex lift of a surface. 
This is an analogue of the twistor lifts and 
a section of the complexification of the two-fold exterior power 
of the pull-back bundle. 
In this paper, 
we will characterize a surface in $N$ such that 
the covariant derivative of the complex lift 
by $\partial /\partial w$ or $\partial /\partial \overline{w}$ 
for every local complex coordinate $w$ of the surface vanishes 
(Theorem~\ref{thm:csect}). 
The twistor lifts of such a surface are degenerate. 
If the mean curvature vector of such a surface vanishes, 
then the surface is totally geodesic. 

Suppose that $N$ is neutral. 
Then analogous discussions and results are valid 
for space-like surfaces in $N$ (see Section~\ref{sect:ns}). 
Refer to \cite{ando4}, \cite{ando5}, \cite{BDM2} 
for space-like twistor spaces. 
For time-like surfaces in $N$, 
we can obtain analogues of Propositions~\ref{pro:WX+YZ=0}, 
\ref{pro:nondegL=0}, \ref{pro:nondegLnot=0}. 
Refer to \cite{HM}, \cite{JR} for time-like twistor spaces. 
In addition, 
we will characterize a time-like surface in $N$ with $L_0 \not= 0$ 
which has local sections of the time-like twistor spaces 
on a neighborhood of each point of the surface 
with fully light-like covariant derivative 
(Theorem~\ref{thm:tnfl}). 
We can refer to \cite{ando8} for such sections. 
The time-like twistor lifts of such surfaces are degenerate. 
Such a surface as in Theorem~\ref{thm:tnfl} does not necessarily have 
zero mean curvature vector, and 
the paracomplex quartic differential defined on the surface is 
zero or null (Remark~\ref{rem:pcqd}). 
We can refer to \cite{ando4}, \cite{ando5}, \cite{ando9}, \cite{ando8} 
for time-like surfaces with zero mean curvature vector 
such that their time-like twistor lifts have 
zero or light-like covariant derivative 
and obtain characterizations 
of such a surface (\cite{ando5}, \cite{ando9}, \cite{AH}). 
We will characterize time-like surfaces in $N$ 
such that their time-like twistor lifts are horizontal along integral curves 
of a light-like one-dimensional distribution (Theorem~\ref{thm:ldcdz}). 
A time-like surface has this property if and only if 
the surface has zero mean curvature vector 
so that the shape operator of any normal vector field is 
zero or light-like, and therefore its time-like twistor lifts have 
zero or light-like covariant derivative (\cite{ando5}). 
We will define the paracomplex lift of a time-like surface in $N$. 
This is an analogue of the time-like twistor lifts 
and a section of the paracomplexification 
of the two-fold exterior power of the pull-back bundle. 
We will characterize a time-like surface in $N$ 
such that the covariant derivative of the paracomplex lift 
by $\partial /\partial           \check{w}$ 
or $\partial /\partial \overline{\check{w}}$ 
for every local paracomplex coordinate $\check{w}$ of the surface vanishes 
(Theorem~\ref{thm:pcsect}). 
The time-like twistor lifts of such a surface are degenerate, and 
the paracomplex quartic differential defined on the surface is 
zero or null (Remark~\ref{rem:pclQ}). 
If the mean curvature vector of such a surface vanishes, 
then a light-like normal vector field is parallel 
(see Remarks~\ref{rem:ldc}, \ref{rem:+=-}), and therefore 
its time-like twistor lifts have 
zero or light-like covariant derivative (\cite{ando5}). 

Suppose that $N$ is Lorentzian. 
Then analogous discussions and results are also valid 
for space-like or time-like surfaces. 
In this case, 
we need the decomposition of the complexification of 
the two-fold exterior power of the pull-back bundle 
on a space-like or time-like surface. 
The two subbundles in this decomposition correspond to 
the eigenvalues $\pm \sqrt{-1}$ of Hodge's $*$-operator. 
These subbundles have complex rank $3$, and 
suitable local frame fields of each subbundle form 
a principal $SO(3, \mbox{\boldmath{$C$}})$-bundle over the surface. 
We will observe that the Levi-Civita connection of $N$ gives 
connections of the two subbundles (Proposition~\ref{pro:connL}). 
Based on it, 
we will compute the covariant derivatives of sections of 
these subbundles, and in particular, 
computing the curvature tensors, 
we will obtain the expressions of the equations of Gauss, Codazzi and Ricci 
of the surface given in \cite{AH} 
(Propositions~\ref{pro:GCRLs}, \ref{pro:GCRLt}). 
The complex twistor spaces associated with the pull-back bundle are 
fiber bundles over the surface in the two subbundles. 
The complex twistor lifts of the surface are 
sections of the complex twistor spaces. 
Then the degeneracy of the complex twistor lifts just means that 
the surface satisfies (i), (ii) in the second paragraph 
of the present section. 
We will obtain characterizations of space-like surfaces 
such that the complex twistor lifts are nondegenerate 
(Propositions~\ref{pro:nondegLsL=0}, \ref{pro:nondegLsLnot=0}). 
In addition, 
we will characterize a space-like surface in $N$ 
such that the covariant derivative of a complex twistor lift 
by $\partial /\partial w$ or $\partial /\partial \overline{w}$ 
for every local complex coordinate $w$ of the surface vanishes 
(Theorem~\ref{thm:deldelbar}). 
Theorem~\ref{thm:deldelbar} seems to be an analogue 
of Theorem~\ref{thm:pcsect} (see Remark~\ref{rem:deldelbar}). 
The complex twistor lifts of such a surface are degenerate, and 
the complex quartic differential defined on the surface is zero 
(Remark~\ref{rem:cqd}). 
We will find a characterization of a space-like surface 
with the above property and zero mean curvature vector 
(see Remark~\ref{rem:ldcLs}, \cite{ando3}, \cite{ando6}, \cite{AH} 
for characterizations of such a surface), 
and directly show that if such a surface has no totally geodesic points, 
then the surface is strictly isotropic, 
which is already obtained in \cite{ando6} (Corollary~\ref{cor:delbar}). 
For a time-like surface in $N$, 
we can obtain analogous results 
related to degeneracy or nondegeneracy of the complex twistor lifts. 
In addition, 
we will characterize a time-like surface in $N$ 
such that a complex twistor lift is horizontal 
along integral curves of a light-like one-dimensional distribution 
(Theorem~\ref{thm:ldcdztL}). 
Theorem~\ref{thm:ldcdztL} is an analogue of Theorem~\ref{thm:ldcdz}. 
A time-like surface has this property if and only if 
the surface has zero mean curvature vector 
so that the shape operator of any normal vector field is 
zero or light-like. 
Referring to \cite{ando6}, 
we can obtain characterizations of such a surface (\cite{AH}). 
The complex twistor lifts of such a surface are degenerate.  

\section{Surfaces 
in Riemannian 4-dimensional space forms}\label{sect:Riem} 

\setcounter{equation}{0} 

\subsection{The equations of Gauss, Codazzi and Ricci} 

Let $N$ be a 4-dimensional Riemannian space form 
with constant sectional curvature $L_0$. 
If $L_0 =0$, 
then we can suppose $N=E^4$; 
if $L_0 >0$, 
then we can suppose $N=\{ x\in E^5 \ | \ \langle x, x\rangle =1/L_0 \}$; 
if $L_0 <0$, 
then we can suppose $N=\{ x\in E^5_1 \ | \ \langle x, x\rangle =1/L_0 \}$, 
where $\langle \ , \ \rangle$ is the metric of $E^5$ or $E^5_1$.  
Let $h$ be the metric of $N$. 
Let $M$ be a Riemann surface 
and $F:M\longrightarrow N$ a conformal immersion of $M$ into $N$. 
Let $(u, v)$ be local isothermal coordinates of $M$ 
compatible with the complex structure of $M$. 
Then the induced metric $g$ on $M$ by $F$ is represented 
as $g=e^{2\lambda} (du^2 +dv^2 )$ 
for a real-valued function $\lambda$. 
Let $\Tilde{\nabla}$ denote the Levi-Civita connection 
of $E^4$, $E^5$ or $E^5_1$ according to $L_0 =0$, $>0$ or $<0$. 
We set $T_1 :=dF(\partial /\partial u )$, 
       $T_2 :=dF(\partial /\partial v )$. 
Let $N_1$, $N_2$ be normal vector fields of $F$ satisfying 
\begin{equation} 
h(N_1 , N_1 )=h(N_2 , N_2 )=e^{2\lambda} , \quad 
h(N_1 , N_2 )=0. 
\label{N1N2} 
\end{equation} 
Then we have 
\begin{equation} 
\begin{split} 
\Tilde{\nabla}_{T_1} (T_1 \ T_2 \ N_1 \ N_2 \ F) 
& =                  (T_1 \ T_2 \ N_1 \ N_2 \ F)S, \\ 
\Tilde{\nabla}_{T_2} (T_1 \ T_2 \ N_1 \ N_2 \ F) 
& =                  (T_1 \ T_2 \ N_1 \ N_2 \ F)T, 
\end{split} 
\label{dt1t2} 
\end{equation} 
where 
\begin{equation} 
\begin{split} 
S & =\left[ 
     \begin{array}{ccccc} 
      \lambda_u         & \lambda_v & -\alpha_1  & -\beta_1   & 1 \\ 
     -\lambda_v         & \lambda_u & -\alpha_2  & -\beta_2   & 0 \\ 
      \alpha_1          & \alpha_2  &  \lambda_u & -\mu_1     & 0 \\ 
      \beta_1           & \beta_2   &  \mu_1     &  \lambda_u & 0 \\ 
      -L_0 e^{2\lambda} &  0        &   0        &   0        & 0 
       \end{array} 
     \right] , \\ 
T & =\left[ 
     \begin{array}{ccccc} 
      \lambda_v & -\lambda_u        & -\alpha_2  & -\beta_2   & 0 \\ 
      \lambda_u &  \lambda_v        & -\alpha_3  & -\beta_3   & 1 \\ 
      \alpha_2  &  \alpha_3         &  \lambda_v & -\mu_2     & 0 \\ 
      \beta_2   &  \beta_3          &  \mu_2     &  \lambda_v & 0 \\ 
       0        & -L_0 e^{2\lambda} &   0        &   0        & 0 
       \end{array} 
    \right] , 
\end{split} 
\label{ST} 
\end{equation} 
and $\alpha_k$, $\beta_k$ ($k=1, 2, 3$) and $\mu_l$ ($l=1, 2$) are 
real-valued functions. 
From \eqref{dt1t2}, 
we obtain $S_v -T_u =ST-TS$ with \eqref{ST}. 
This is equivalent to 
the system of the equations of Gauss, Codazzi and Ricci. 
The equation of Gauss is given by 
\begin{equation} 
  \lambda_{uu} +\lambda_{vv} +L_0 e^{2\lambda} 
=-\alpha_1 \alpha_3 -\beta_1 \beta_3 +\alpha^2_2 +\beta^2_2 . 
\label{gauss}
\end{equation}
The equations of Codazzi are given by 
\begin{equation} 
\begin{split}
  (\alpha_1 )_v -(\alpha_2 )_u 
& =\alpha_2 \lambda_u +\alpha_3 \lambda_v 
   -\beta_2 \mu_1     +\beta_1  \mu_2  , \\ 
   (\alpha_2 )_v -(\alpha_3 )_u 
& =-\alpha_1 \lambda_u -\alpha_2 \lambda_v 
   -\beta_3  \mu_1     +\beta_2  \mu_2  , \\ 
   (\beta_1 )_v -(\beta_2 )_u 
& = \beta_2  \lambda_u +\beta_3  \lambda_v 
   +\alpha_2 \mu_1     -\alpha_1 \mu_2  , \\ 
   (\beta_2 )_v -(\beta_3 )_u 
& =-\beta_1  \lambda_u -\beta_2  \lambda_v 
   +\alpha_3 \mu_1     -\alpha_2 \mu_2 .  
\end{split} 
\label{codazzi} 
\end{equation} 
The equation of Ricci is given by 
\begin{equation} 
  (\mu_1 )_v -(\mu_2 )_u = \alpha_1 \beta_2 -\alpha_2 \beta_1 
                          +\alpha_2 \beta_3 -\alpha_3 \beta_2 . 
\label{ricci}
\end{equation}

Suppose that $N$ is oriented. 
The metric $h$ of $N$ induces a metric $\hat{h}$ 
of the two-fold exterior power $\bigwedge^2\!F^*T\!N$ 
of the pull-back bundle $F^*T\!N$ on $M$ by $F$. 
The bundle $\bigwedge^2\!F^*T\!N$ is of rank $6$ and decomposed into 
two subbundles $\bigwedge^2_{\pm}\!F^*T\!N$ of rank $3$. 
These subbundles are orthogonal to each other with respect to $\hat{h}$ 
and correspond to the eigenvalues $\pm 1$ respectively 
of Hodge's $*$-operator with respect to $h$ and the orientation of $N$. 
Suppose that $(T_1 , T_2 , N_1 , N_2 )$ gives the orientation of $N$. 
We set 
\begin{equation} 
e_1 :=\dfrac{1}{e^{\lambda}} T_1 , \quad 
e_2 :=\dfrac{1}{e^{\lambda}} T_2 , \quad 
e_3 :=\dfrac{1}{e^{\lambda}} N_1 , \quad 
e_4 :=\dfrac{1}{e^{\lambda}} N_2 
\label{e1e2e3e4} 
\end{equation} 
and 
\begin{equation} 
\begin{split} 
    \Theta_{\pm , 1} 
& :=\dfrac{1}{\sqrt{2}} (e_1 \wedge e_2 \pm e_3 \wedge e_4 ), \\ 
    \Theta_{\pm , 2} 
& :=\dfrac{1}{\sqrt{2}} (e_1 \wedge e_3 \pm e_4 \wedge e_2 ), \\ 
    \Theta_{\pm , 3} 
& :=\dfrac{1}{\sqrt{2}} (e_1 \wedge e_4 \pm e_2 \wedge e_3 ). 
\end{split} 
\label{Omega} 
\end{equation} 
Then $\Theta_{\pm , 1}$, $\Theta_{\pm , 2}$, $\Theta_{\pm , 3}$ form 
orthonormal local frame fields of $\bigwedge^2_{\pm}\!F^*T\!N$ 
respectively. 
The \textit{twistor lifts\/} of $F$ are 
sections of the twistor spaces associated with 
the pull-back bundle $F^*T\!N$ and locally given by $\Theta_{\pm , 1}$, 
which are determined by the immersion $F$ 
and depend on the choices of neither the isothermal coordinates $(u, v)$ 
nor normal vector fields $N_1$, $N_2$. 

Let $\nabla$ be the Levi-Civita connection of $h$ 
and $\hat{\nabla}$ the induced connection of $\bigwedge^2\!F^*T\!N$ 
by $\nabla$. 
Then $\hat{\nabla}$ gives connections of $\bigwedge^2_{\pm}\!F^*T\!N$. 
By \eqref{dt1t2} with \eqref{ST}, we obtain 
\begin{equation} 
\begin{split} 
    \hat{\nabla}_{T_1}  
   (\Theta_{\pm , 1} \ \Theta_{\pm , 2} \ \Theta_{\pm , 3} ) 
& =(\Theta_{\pm , 1} \ \Theta_{\pm , 2} \ \Theta_{\pm , 3} ) 
    \left[ \begin{array}{ccc} 
            0       &  -   W_{\pm}   &  -   Y_{\mp}   \\ 
            W_{\pm} &      0         & \pm \psi_{\pm} \\ 
            Y_{\mp} & \mp \psi_{\pm} &      0 
             \end{array} 
    \right] , \\ 
    \hat{\nabla}_{T_2} 
   (\Theta_{\pm , 1} \ \Theta_{\pm , 2} \ \Theta_{\pm , 3} ) 
& =(\Theta_{\pm , 1} \ \Theta_{\pm , 2} \ \Theta_{\pm , 3} ) 
    \left[ \begin{array}{ccc} 
            0       & \mp  Z_{\pm}   & \pm  X_{\mp}   \\ 
        \pm Z_{\pm} &      0         & \mp \phi_{\mp} \\ 
        \mp X_{\mp} & \pm \phi_{\mp} &      0 
             \end{array} 
    \right] , 
\end{split} 
\label{hatnabla} 
\end{equation} 
where 
\begin{equation} 
\begin{array}{lcl} 
W_{\pm} :=\alpha_2 \pm \beta_1  , & \ & X_{\pm} :=\alpha_2 \pm \beta_3 , \\ 
Y_{\pm} :=\beta_2  \pm \alpha_1 , & \ & Z_{\pm} :=\beta_2  \pm \alpha_3 
\end{array} 
\label{WXYZ} 
\end{equation} 
and 
\begin{equation} 
\phi_{\pm} :=\lambda_u \mp \mu_2 , \quad 
\psi_{\pm} :=\lambda_v \mp \mu_1 . 
\label{phipsi} 
\end{equation} 
We have 
\begin{equation} 
W_+ +W_- =X_+ +X_- , \quad 
Y_+ +Y_- =Z_+ +Z_- . 
\label{W=XY=Z} 
\end{equation} 

Let $\hat{R}$ be the curvature tensor of $\hat{\nabla}$. 
Since $N$ is a space form of constant sectional curvature $L_0$, 
we obtain 
\begin{equation} 
  \hat{R} (T_1 , T_2 ) 
 (\Theta_{\pm , 1} \ \Theta_{\pm , 2} \ \Theta_{\pm , 3} ) 
=(\Theta_{\pm , 1} \ \Theta_{\pm , 2} \ \Theta_{\pm , 3} ) 
  \left[ \begin{array}{ccc} 
          0 &     0                &     0                \\ 
          0 &     0                & \pm L_0 e^{2\lambda} \\ 
          0 & \mp L_0 e^{2\lambda} &     0 
           \end{array} 
  \right] . 
\label{hatR} 
\end{equation} 
On the other hand, 
we can compute the left side of \eqref{hatR} 
by \eqref{hatnabla}, and then we obtain  

\begin{pro}\label{pro:GCRR} 
The functions $W_{\pm}$, $X_{\pm}$, $Y_{\pm}$, $Z_{\pm}$ 
as in \eqref{WXYZ} satisfy not only \eqref{W=XY=Z} but also 
\begin{equation} 
 W_{\mp} X_{\pm} +Y_{\pm} Z_{\mp} 
=L_0 e^{2\lambda} +(\phi_{\pm} )_u +(\psi_{\mp} )_v 
\label{gaussricci} 
\end{equation} 
and 
\begin{equation} 
\begin{split}
      (Y_{\pm} )_v \mp (X_{\pm} )_u 
& =\pm W_{\mp} \phi_{\pm} -Z_{\mp} \psi_{\mp} , \\  
      (W_{\mp} )_v \pm (Z_{\mp} )_u 
& =\mp Y_{\pm} \phi_{\pm} -X_{\pm} \psi_{\mp} 
\end{split} 
\label{codazziWXYZ} 
\end{equation} 
with \eqref{phipsi}. 
\end{pro} 

\begin{rem} 
By direct computations as in \cite{AH}, 
we see that 
the system of two equations in \eqref{gaussricci} is equivalent to 
the system of equations \eqref{gauss}, \eqref{ricci} and that 
the system of four equations in \eqref{codazziWXYZ} is equivalent to 
the system of four equations in \eqref{codazzi}. 
\end{rem} 

\subsection{Degeneracy of the twistor lifts} 

Let $K$ be the curvature of $g$. 
Then we have $K=-e^{-2\lambda} (\lambda_{uu} +\lambda_{vv} )$. 
Let $\nabla^{\perp}$ be the normal connection of $F$ 
and $R^{\perp}$ the curvature tensor of $\nabla^{\perp}$. 
Then $R^{\perp} \equiv 0$ just means that $\nabla^{\perp}$ is flat, 
which is equivalent to $(\mu_1 )_v =(\mu_2 )_u$. 
The twistor lifts $\Theta_{\pm , 1}$ of $F$ are said to be 
\textit{nondegenerate\/} (respectively, \textit{degenerate\/}) 
if for each $\varepsilon \in \{ +, -\}$, 
    $\hat{\nabla}_{T_1} \Theta_{\varepsilon , 1}$ 
and $\hat{\nabla}_{T_2} \Theta_{\varepsilon , 1}$ are 
linearly independent (respectively, dependent) at each point of $M$. 
From \eqref{hatnabla} and \eqref{gaussricci}, 
we obtain 

\begin{pro}\label{pro:WX+YZ=0} 
The immersion $F$ satisfies both $K\equiv L_0$ and $R^{\perp} \equiv 0$ 
if and only if the twistor lifts $\Theta_{\pm , 1}$ of $F$ are 
degenerate. 
\end{pro} 

The twistor lifts $\Theta_{\pm , 1}$ of $F$ are nondegenerate 
if and only if $\Delta_{\pm} :=W_{\mp} X_{\pm} +Y_{\pm} Z_{\mp} \not= 0$. 
If we suppose $\Delta_{\pm} \not= 0$, 
then \eqref{codazziWXYZ} can be rewritten into 
$(\phi_{\pm} , \psi_{\mp} )=(A_{\pm} , B_{\mp} )$, 
where 
\begin{equation} 
  \left[ \begin{array}{c} 
          A_{\pm} \\ 
          B_{\mp} 
           \end{array} 
  \right] 
:=\mp 
  \dfrac{1}{\Delta_{\pm}} 
  \left[ \begin{array}{cc} 
         -  X_{\pm} &     Z_{\mp} \\ 
        \pm Y_{\pm} & \pm W_{\mp} 
           \end{array} 
  \right] 
  \left[ \begin{array}{c} 
          (Y_{\pm} )_v \mp (X_{\pm} )_u \\ 
          (W_{\mp} )_v \pm (Z_{\mp} )_u 
           \end{array} 
  \right] . 
\label{ABWXYZ} 
\end{equation} 

Let $W_{\pm}$, $X_{\pm}$, $Y_{\pm}$, $Z_{\pm}$ be functions 
of two variables $u$, $v$ satisfying \eqref{W=XY=Z}, 
$\Delta_{\pm} \not= 0$, 
\begin{equation} 
(A_{\pm} )_u +(B_{\mp} )_v =\Delta_{\pm} 
\label{gaussricciRAB} 
\end{equation} 
and 
\begin{equation} 
(A_+ +A_- )_v =(B_+ +B_- )_u , 
\label{lambdauv} 
\end{equation} 
where $A_{\pm}$, $B_{\mp}$ are functions 
constructed by $W_{\mp}$, $X_{\pm}$, $Y_{\pm}$, $Z_{\mp}$ 
as in \eqref{ABWXYZ}. 
Then $W_{\pm}$, $X_{\pm}$, $Y_{\pm}$, $Z_{\pm}$, 
$\phi_{\pm} :=A_{\pm}$, $\psi_{\mp} :=B_{\mp}$ satisfy 
\eqref{gaussricci} with $L_0 =0$ and \eqref{codazziWXYZ}. 
By \eqref{lambdauv}, 
there exists a function $\lambda$ satisfying 
\begin{equation} 
\lambda_u =\dfrac{1}{2} (A_+ +A_- ), \quad 
\lambda_v =\dfrac{1}{2} (B_+ +B_- ). 
\label{lambdaAB} 
\end{equation} 
We set 
\begin{equation} 
\begin{array}{lclcl} 
\alpha_1 :=\dfrac{1}{2} (Y_+ -Y_- ), & \ & 
\alpha_2 :=\dfrac{1}{2} (X_+ +X_- ), & \ & 
\alpha_3 :=\dfrac{1}{2} (Z_+ -Z_- ), \\ 
\ & \ & \ & \ & \ \\ 
\beta_1  :=\dfrac{1}{2} (W_+ -W_- ), & \ & 
\beta_2  :=\dfrac{1}{2} (Y_+ +Y_- ), & \ & 
\beta_3  :=\dfrac{1}{2} (X_+ -X_- ) 
\end{array} 
\label{alphabeta} 
\end{equation} 
and 
\begin{equation} 
\mu_1 :=\dfrac{1}{2} (B_- -B_+ ), \quad 
\mu_2 :=\dfrac{1}{2} (A_- -A_+ ). 
\label{mu1mu2} 
\end{equation} 
Then $\lambda$, $\alpha_k$, $\beta_k$ ($k=1, 2, 3$), $\mu_l$ ($l=1, 2$) 
satisfy \eqref{gauss} with $L_0 =0$, \eqref{codazzi}, \eqref{ricci}. 
Therefore these functions determine an immersion $F$ 
into $E^4$ such that the twistor lifts of $F$ are nondegenerate, 
which is unique up to an isometry of $E^4$. 

Hence we obtain 

\begin{pro}\label{pro:nondegL=0} 
If the twistor lifts $\Theta_{\pm , 1}$ of the immersion $F$ 
into $E^4$ are nondegenerate, 
then $\Delta_{\pm} \not= 0$, 
and $A_{\pm}$, $B_{\mp}$ as in \eqref{ABWXYZ} satisfy 
$(A_{\pm} , B_{\mp} )=(\phi_{\pm} , \psi_{\mp} )$ with \eqref{phipsi} 
and \eqref{gaussricciRAB}. 
In addition, 
functions $W_{\pm}$, $X_{\pm}$, $Y_{\pm}$, $Z_{\pm}$ 
of two variables $u$, $v$ satisfying \eqref{W=XY=Z}, 
$\Delta_{\pm} \not= 0$, \eqref{gaussricciRAB} and \eqref{lambdauv} 
give an immersion $F$ into $E^4$ 
such that the twistor lifts are nondegenerate, 
which is unique up to an isometry of $E^4$. 
\end{pro} 

Let $W_{\pm}$, $X_{\pm}$, $Y_{\pm}$, $Z_{\pm}$ be functions 
of two variables satisfying \eqref{W=XY=Z}, 
$\Delta_{\pm} \not= 0$ and 
\begin{equation} 
\Delta_+ -\Delta_- =(A_+ -A_- )_u -(B_+ -B_- )_v , 
\label{Delta+Delta-} 
\end{equation} 
and suppose that 
\begin{equation} 
f:=\Delta_+ -(A_+ )_u -(B_- )_v \ (=\Delta_- -(A_- )_u -(B_+ )_v ) 
\label{f} 
\end{equation} 
is nowhere zero and satisfies 
\begin{equation} 
f_u =f(A_+ +A_- ), \quad 
f_v =f(B_+ +B_- ). 
\label{fufv} 
\end{equation} 
For a nonzero number $L_0$, 
suppose $f/L_0 >0$ and set $\lambda :=(1/2)\log f/L_0$. 
Then $\lambda$ satisfies \eqref{lambdaAB}. 
Let $\mu_1$, $\mu_2$ be functions as in \eqref{mu1mu2}. 
Then $\phi_{\pm}$, $\psi_{\pm}$ as in \eqref{phipsi} 
satisfy $\phi_{\pm} =A_{\pm}$, $\psi_{\mp} =B_{\mp}$. 
Therefore we obtain \eqref{gaussricci} and \eqref{codazziWXYZ}, 
and for functions $\alpha_k$, $\beta_k$ ($k=1, 2, 3$) 
given as in \eqref{alphabeta}, 
we obtain \eqref{gauss}, \eqref{codazzi}, \eqref{ricci}. 
Therefore the functions 
$\lambda$, $\alpha_k$, $\beta_k$ ($k=1, 2, 3$), $\mu_l$ ($l=1, 2$) 
determine an immersion $F$ into 
a 4-dimensional Riemannian space form $N$ 
with constant sectional curvature $L_0$ 
such that the twistor lifts of $F$ are nondegenerate, 
and this immersion is unique up to an isometry of $N$. 

Hence we obtain 

\begin{pro}\label{pro:nondegLnot=0} 
If the twistor lifts $\Theta_{\pm , 1}$ of the immersion $F$ 
into $N$ with $L_0 \not= 0$ are nondegenerate, 
then $\Delta_{\pm} \not= 0$, 
and $A_{\pm}$, $B_{\mp}$ as in \eqref{ABWXYZ} satisfy 
$(A_{\pm} , B_{\mp} )=(\phi_{\pm} , \psi_{\mp} )$ with \eqref{phipsi} 
and \eqref{fufv} with $f=L_0 e^{2\lambda}$. 
In addition, 
a nonzero number $L_0$ and 
functions $W_{\pm}$, $X_{\pm}$, $Y_{\pm}$, $Z_{\pm}$ 
of two variables $u$, $v$ satisfying \eqref{W=XY=Z}, 
$\Delta_{\pm} \not= 0$, \eqref{Delta+Delta-} 
and \eqref{fufv} with \eqref{f} and $f/L_0 >0$ 
give an immersion $F$ into $N$ with constant sectional curvature $L_0$ 
such that the twistor lifts are nondegenerate, 
which is unique up to an isometry of $N$. 
\end{pro} 

\subsection{Horizontality of the complex lift} 

Let $\Theta$ be a section of 
the complexification $\bigwedge^2\!F^*T\!N \otimes \mbox{\boldmath{$C$}}$ 
of $\bigwedge^2\!F^*T\!N$ which is locally given by 
\begin{equation} 
\begin{split} 
     \Theta 
:= & \dfrac{1}{\sqrt{2}} (e_1 \wedge e_2 +\sqrt{-1}  e_3 \wedge e_4 ) \\ 
 = & \dfrac{1}{2} ( \Theta_{+, 1} +\Theta_{-, 1} 
     +\sqrt{-1}    (\Theta_{+, 1} -\Theta_{-, 1} )).  
\end{split} 
\label{ThetaR} 
\end{equation} 
Then $\Theta$ in \eqref{ThetaR} is determined by the immersion $F$ 
and the orientations of $M$, $N$. 
We call $\Theta$ the \textit{complex lift\/} of $F$. 
Let $w:=u+\sqrt{-1} v$ be a local complex coordinate of $M$ 
and set $\partial 
       :=\partial /\partial w 
        =(1/2)(\partial /\partial u -\sqrt{-1} \partial /\partial v)$. 
Then we can naturally consider 
the covariant derivative $\hat{\nabla}_{\partial} \Theta$ 
of $\Theta$ by $\partial$ and we have 
\begin{equation} 
\begin{split} 
    \hat{\nabla}_{\partial} \Theta 
= & \dfrac{1}{4} (  \hat{\nabla}_{T_1} (\Theta_{+, 1} +\Theta_{-, 1} ) 
                   +\hat{\nabla}_{T_2} (\Theta_{+, 1} -\Theta_{-, 1} ) \\ 
  & \quad          +\sqrt{-1} 
                  ( \hat{\nabla}_{T_1} (\Theta_{+, 1} -\Theta_{-, 1} ) 
                   -\hat{\nabla}_{T_2} (\Theta_{+, 1} +\Theta_{-, 1} ))).  
\end{split} 
\label{hndT} 
\end{equation} 
Whether $\hat{\nabla}_{\partial} \Theta$ vanishes does not depend on 
the choice of $w$, and 
from \eqref{hndT}, we see that $\hat{\nabla}_{\partial} \Theta$ vanishes 
if and only if the following hold: 
\begin{equation} 
  \hat{\nabla}_{T_1} \Theta_{+, 1} 
= \hat{\nabla}_{T_2} \Theta_{-, 1} , \quad 
  \hat{\nabla}_{T_2} \Theta_{+, 1} 
=-\hat{\nabla}_{T_1} \Theta_{-, 1} . 
\label{hndT0} 
\end{equation} 
By \eqref{hatnabla}, \eqref{hndT0} is equivalent to 
\begin{equation} 
(W_{\pm} , X_{\pm} )=(-Z_{\mp} , Y_{\mp} ). 
\label{WX-ZY} 
\end{equation} 
Suppose \eqref{WX-ZY}. 
Then by \eqref{W=XY=Z}, we have 
\begin{equation} 
W_+ =-W_- = Z_+ =-Z_- , \quad 
X_+ =-X_- =-Y_+ = Y_- . 
\label{W+-Z+-X+-Y+-} 
\end{equation} 
Applying \eqref{W+-Z+-X+-Y+-} to \eqref{gaussricci}, 
we obtain $K\equiv L_0$ and $R^{\perp} \equiv 0$. 
The former condition is equivalent to 
\begin{equation} 
\lambda_{uu} +\lambda_{vv} +L_0 e^{2\lambda} =0 
\label{K=L0} 
\end{equation} 
and the latter condition is equivalent to the existence of 
a function $\gamma$ satisfying 
\begin{equation} 
\gamma_u =\mu_1 , \quad 
\gamma_v =\mu_2 . 
\label{gamma} 
\end{equation} 
Applying \eqref{W+-Z+-X+-Y+-} to \eqref{codazziWXYZ}, we obtain 
\begin{equation} 
  (\partial_u \pm \partial_v )X_+ 
= (\phi_{\pm} \mp \psi_{\mp} )W_+ , \quad 
  (\partial_u \pm \partial_v )W_+ 
=-(\phi_{\pm} \mp \psi_{\mp} )X_+ , 
\label{dupmdv} 
\end{equation} 
where $\partial_u :=\partial /\partial u$, 
      $\partial_v :=\partial /\partial v$. 
Applying \eqref{dupmdv} to 
\begin{equation*} 
 (\partial_u +\partial_v )(\partial_u -\partial_v )X_+ 
=(\partial_u -\partial_v )(\partial_u +\partial_v )X_+ , 
\end{equation*} 
we obtain 
\begin{equation} 
 W_+ (\partial_u +\partial_v )(\phi_- +\psi_+ ) 
=W_+ (\partial_u -\partial_v )(\phi_+ -\psi_- ); 
\label{W+} 
\end{equation} 
applying \eqref{dupmdv} to 
\begin{equation*} 
 (\partial_u +\partial_v )(\partial_u -\partial_v )W_+ 
=(\partial_u -\partial_v )(\partial_u +\partial_v )W_+ , 
\end{equation*} 
we obtain 
\begin{equation} 
 X_+ (\partial_u +\partial_v )(\phi_- +\psi_+ ) 
=X_+ (\partial_u -\partial_v )(\phi_+ -\psi_- ). 
\label{X+} 
\end{equation} 
If $F$ satisfies 
\begin{equation} 
 (\partial_u +\partial_v )(\phi_- +\psi_+ ) 
=(\partial_u -\partial_v )(\phi_+ -\psi_- ), 
\label{dupmdv2} 
\end{equation} 
then by \eqref{phipsi} and \eqref{gamma}, 
we obtain $\lambda_{uv} =0$. 
Therefore by \eqref{K=L0}, 
if $L_0 =0$, 
then $\lambda$ is given by $\lambda =au+bv+c$ for constants $a$, $b$, $c$; 
if $L_0 \not= 0$, 
then $\lambda$ is a function of one variable $u$ or $v$: 
$\lambda =p(u)$ or $q(v)$, 
where $p$, $q$ are solutions of 
an ordinary differential equation $f'' +L_0 e^{2f} =0$. 

Let $\lambda$ be a function given by $au+bv+c$, $p(u)$ or $q(v)$. 
Let $\gamma$  be a function of two variables $u$, $v$. 
Then \eqref{dupmdv} with \eqref{phipsi} and \eqref{gamma} has 
a unique solution $(W_+ , X_+ )$ 
for an arbitrarily given initial value. 
Let $W_-$, $X_-$, $Y_{\pm}$, $Z_{\pm}$ be functions of $u$, $v$ 
given by \eqref{W+-Z+-X+-Y+-}. 
Then these functions satisfy \eqref{W=XY=Z}, \eqref{gaussricci} 
and \eqref{codazziWXYZ}. 
Therefore these functions give an immersion into $N$ 
satisfying $\hat{\nabla}_{\partial} \Theta =0$ 
for $\Theta$ in \eqref{ThetaR}, 
which is unique up to an isometry of $N$. 

We can similarly study 
the condition $\hat{\nabla}_{\overline{\partial}} \Theta =0$  
for $\overline{\partial} 
   :=\partial /\partial \overline{w} 
    =(1/2)(\partial /\partial u +\sqrt{-1} \partial /\partial v)$. 
This condition is equivalent to 
\begin{equation} 
  \hat{\nabla}_{T_1} \Theta_{+, 1} 
=-\hat{\nabla}_{T_2} \Theta_{-, 1} , \quad 
  \hat{\nabla}_{T_2} \Theta_{+, 1} 
= \hat{\nabla}_{T_1} \Theta_{-, 1} , 
\label{hndT0bar} 
\end{equation} 
instead of \eqref{hndT0}. 
Therefore \eqref{hndT0bar} is equivalent to 
\begin{equation} 
(W_{\pm} , X_{\pm} )=(Z_{\mp} , -Y_{\mp} ), 
\label{WX-ZYbar} 
\end{equation} 
instead of \eqref{WX-ZY}. 
If we suppose \eqref{WX-ZYbar}, 
then by \eqref{W=XY=Z}, we have 
\begin{equation} 
W_+ =-W_- =-Z_+ = Z_- , \quad 
X_+ =-X_- = Y_+ =-Y_- , 
\label{W+-Z+-X+-Y+-bar} 
\end{equation} 
instead of \eqref{W+-Z+-X+-Y+-}, 
and applying \eqref{W+-Z+-X+-Y+-bar} to \eqref{codazziWXYZ}, we obtain 
\begin{equation} 
  (\partial_u \mp \partial_v )X_+ 
= (\phi_{\pm} \pm \psi_{\mp} )W_+ , \quad 
  (\partial_u \mp \partial_v )W_+ 
=-(\phi_{\pm} \pm \psi_{\mp} )X_+ , 
\label{dupmdvbar} 
\end{equation} 
instead of \eqref{dupmdv}. 

Hence we obtain 

\begin{thm}\label{thm:csect} 
The complex lift $\Theta$ defined in \eqref{ThetaR} satisfies 
either $\hat{\nabla}_{\partial} \Theta =0$ 
or     $\hat{\nabla}_{\overline{\partial}} \Theta =0$ 
if and only if the immersion $F$ satisfies 
\eqref{W+-Z+-X+-Y+-} or \eqref{W+-Z+-X+-Y+-bar}. 
If $F$ satisfies \eqref{W+-Z+-X+-Y+-} or \eqref{W+-Z+-X+-Y+-bar}, 
then 
\begin{itemize} 
\item[{\rm (a)}]{$K\equiv L_0$ and $R^{\perp} \equiv 0$,} 
\item[{\rm (b)}]{$W_+$, $X_+$ satisfy 
either \eqref{dupmdv} or \eqref{dupmdvbar}.} 
\end{itemize} 
In addition, 
for a function $\lambda =au+bv+c$, $p(u)$ or $q(v)$, 
a function $\gamma$ of two variables $u$, $v$ and 
a unique solution $(W_+ , X_+ )$ of \eqref{dupmdv} or \eqref{dupmdvbar} 
determined by an arbitrarily given initial value, 
an immersion into $N$ satisfying 
either $\hat{\nabla}_{\partial} \Theta =0$ 
or     $\hat{\nabla}_{\overline{\partial}} \Theta =0$ is constructed 
by \eqref{W+-Z+-X+-Y+-} or \eqref{W+-Z+-X+-Y+-bar}, 
and unique up to an isometry of $N$. 
\end{thm} 

\begin{rem} 
By \eqref{WXYZ}, 
the mean curvature vector of $F$ vanishes 
if and only if $F$ satisfies 
\begin{equation} 
(W_{\pm} , Z_{\pm} )=(X_{\mp} , Y_{\mp} ), 
\label{WZXY} 
\end{equation} 
and if we suppose \eqref{WZXY} and one of \eqref{WX-ZY}, \eqref{WX-ZYbar}, 
then we obtain $W_{\pm} =X_{\pm} =Y_{\pm} =Z_{\pm} =0$. 
\end{rem} 

\begin{rem} 
We say that the second fundamental form $\sigma$ of $F$ satisfies 
the \textit{linearly dependent condition\/} 
if $F$ satisfies 
\begin{equation} 
 \cos \theta (\alpha_1 , \alpha_2 , \alpha_3 ) 
+\sin \theta (\beta_1  , \beta_2  , \beta_3  ) 
=0
\label{ldc} 
\end{equation} 
for a function $\theta$ of $u$, $v$. 
If functions $W_+$, $X_+$ are nonzero, 
then the second fundamental form of an immersion into $N$ 
as in Theorem~\ref{thm:csect} does not satisfy 
the linearly dependent condition and therefore 
the image is not contained 
in any totally geodesic hypersurface of $N$. 
Refer to \cite{AH} for the linearly dependent condition. 
\end{rem} 

\begin{rem}\label{rem:cqdR} 
Let $F:M\longrightarrow N$ be a conformal immersion. 
Let $w=u+\sqrt{-1} v$ be a local complex coordinate of $M$. 
Let $Q$ be a complex quartic differential on $M$ defined by 
\begin{equation} 
Q=h(\sigma (\partial , \partial ), 
    \sigma (\partial , \partial ))dw^4 , 
\label{QR} 
\end{equation} 
where $\sigma$ is the second fundamental form of $F$. 
Then we obtain 
\begin{equation} 
  h(\sigma (\partial , \partial ), 
    \sigma (\partial , \partial )) 
  = \dfrac{e^{2\lambda}}{16} 
    (p -2\sqrt{-1} q)(\overline{p} -2\sqrt{-1} \overline{q} ) 
\label{QRc} 
\end{equation} 
for 
\begin{equation*} 
p:=\alpha_1 -\alpha_3 +\sqrt{-1} (\beta_1 -\beta_3 ), \quad 
q:=\alpha_2           +\sqrt{-1}  \beta_2 . 
\end{equation*} 
Therefore $Q$ is zero if and only if 
one of $p \pm 2\sqrt{-1} q$ is zero. 
Noticing 
\begin{equation*} 
\begin{split} 
    p +2\sqrt{-1} q 
& =-Y_- -Z_+ +\sqrt{-1} (W_+ +X_- ), \\ 
    p -2\sqrt{-1} q 
& = Y_+ +Z_- -\sqrt{-1} (W_- +X_+ ), 
\end{split} 
\end{equation*} 
we see that 
if $F$ satisfies \eqref{W+-Z+-X+-Y+-} or \eqref{W+-Z+-X+-Y+-bar} and 
if $F$ has no totally geodesic points, 
then $Q$ is nowhere zero. 
\end{rem} 
 
\section{Space-like surfaces 
in neutral 4-dimensional space forms}\label{sect:ns} 

\setcounter{equation}{0} 

Let $N$ be a 4-dimensional neutral space form 
with constant sectional curvature $L_0$. 
If $L_0 =0$, 
then we can suppose $N=E^4_2$; 
if $L_0 >0$, 
then we can suppose $N=\{ x\in E^5_2 \ | \ \langle x, x\rangle =1/L_0 \}$; 
if $L_0 <0$, 
then we can suppose $N=\{ x\in E^5_3 \ | \ \langle x, x\rangle =1/L_0 \}$, 
where $\langle \ , \ \rangle$ is the metric of $E^5_2$ or $E^5_3$.  
We denote $N$ with $L_0 =1$, $-1$ by $S^4_2$, $H^4_2$ respectively. 
Let $h$ be the neutral metric of $N$. 
Let $M$ be a Riemann surface 
and $F:M\longrightarrow N$ a space-like and conformal immersion 
of $M$ into $N$. 
Let $(u, v)$, $g$, $\lambda$ be 
as in the beginning of Section~\ref{sect:Riem}. 
Let $\Tilde{\nabla}$ denote the Levi-Civita connection 
of $E^4_2$, $E^5_2$ or $E^5_3$ according to $L_0 =0$, $>0$ or $<0$. 
Let $T_1$, $T_2$ be as in Section~\ref{sect:Riem}. 
Let $N_1$, $N_2$ be normal vector fields of $F$ satisfying 
$$h(N_1 , N_1 )=h(N_2 , N_2 )=-e^{2\lambda} , \quad 
  h(N_1 , N_2 )=0.$$ 
Then we have \eqref{dt1t2} with 
\begin{equation} 
\begin{split} 
S & =\left[ \begin{array}{ccccc} 
         \lambda_u         & \lambda_v &  \alpha_1  &  \beta_1   & 1 \\ 
        -\lambda_v         & \lambda_u &  \alpha_2  &  \beta_2   & 0 \\ 
         \alpha_1          & \alpha_2  &  \lambda_u & -\mu_1     & 0 \\ 
         \beta_1           & \beta_2   &  \mu_1     &  \lambda_u & 0 \\ 
         -L_0 e^{2\lambda} &  0        &   0        &   0        & 0 
              \end{array} 
     \right] , \\ 
T & =\left[ \begin{array}{ccccc} 
         \lambda_v & -\lambda_u        &  \alpha_2  &  \beta_2   & 0 \\ 
         \lambda_u &  \lambda_v        &  \alpha_3  &  \beta_3   & 1 \\ 
         \alpha_2  &  \alpha_3         &  \lambda_v & -\mu_2     & 0 \\ 
         \beta_2   &  \beta_3          &  \mu_2     &  \lambda_v & 0 \\ 
          0        & -L_0 e^{2\lambda} &   0        &   0        & 0 
              \end{array} 
    \right] , 
\end{split} 
\label{STns} 
\end{equation} 
and $\alpha_k$, $\beta_k$ ($k=1, 2, 3$) and $\mu_l$ ($l=1, 2$) are 
real-valued functions. 
From \eqref{dt1t2}, 
we obtain $S_v -T_u =ST-TS$ with \eqref{STns}. 
This is equivalent to 
the system of the equations of Gauss, Codazzi and Ricci. 
The equation of Gauss is given by 
\begin{equation} 
  \lambda_{uu} +\lambda_{vv} +L_0 e^{2\lambda} 
= \alpha_1 \alpha_3 +\beta_1 \beta_3 -\alpha^2_2 -\beta^2_2 . 
\label{gaussns}
\end{equation}
The equations of Codazzi are given by \eqref{codazzi}. 
The equation of Ricci is given by 
\begin{equation} 
  (\mu_1 )_v -(\mu_2 )_u =-\alpha_1 \beta_2 +\alpha_2 \beta_1 
                          -\alpha_2 \beta_3 +\alpha_3 \beta_2 . 
\label{riccins} 
\end{equation}

Suppose that $N$ is oriented. 
Let $\hat{h}$ be the metric of $\bigwedge^2\!F^*T\!N$ induced by $h$. 
Then $\hat{h}$ has signature $(2, 4)$. 
The bundle $\bigwedge^2\!F^*T\!N$ is decomposed into 
two subbundles $\bigwedge^2_{\pm}\!F^*T\!N$ of rank $3$, and 
these subbundles are orthogonal to each other with respect to $\hat{h}$ 
and correspond to the eigenvalues $\pm 1$ respectively 
of Hodge's $*$-operator 
with respect to the neutral metric $h$ and the orientation of $N$. 
The restrictions of $\hat{h}$ on $\bigwedge^2_{\pm}\!F^*T\!N$ have 
signature $(1, 2)$. 
Suppose that $(T_1 , T_2 , N_1 , N_2 )$ gives the orientation of $N$. 
Let $e_k$ ($k=1, 2, 3, 4$) and $\Theta_{\pm , l}$ ($l=1, 2, 3$) be 
as in \eqref{e1e2e3e4}, \eqref{Omega}, respectively. 
Then $\Theta_{\mp , 1}$, $\Theta_{\pm , 2}$, $\Theta_{\pm , 3}$ form 
pseudo-orthonormal local frame fields of $\bigwedge^2_{\pm}\!F^*T\!N$ 
respectively.  
The \textit{space-like twistor lifts\/} of $F$ are 
sections of the space-like twistor spaces associated with $F^*T\!N$, 
and as in Section~\ref{sect:Riem}, 
these sections are locally given by $\Theta_{\mp , 1}$ 
and determined by $F$. 

Let $\nabla$, $\hat{\nabla}$ be as in Section~\ref{sect:Riem}. 
Then $\hat{\nabla}$ gives connections of $\bigwedge^2_{\pm}\!F^*T\!N$. 
By \eqref{dt1t2} with \eqref{STns}, we obtain 
\begin{equation} 
\begin{split} 
    \hat{\nabla}_{T_1}  
   (\Theta_{\mp , 1} \ \Theta_{\pm , 2} \ \Theta_{\pm , 3} ) 
& =(\Theta_{\mp , 1} \ \Theta_{\pm , 2} \ \Theta_{\pm , 3} ) 
    \left[ \begin{array}{ccc} 
            0       &      W_{\pm}   &      Y_{\mp}   \\ 
            W_{\pm} &      0         & \pm \psi_{\pm} \\ 
            Y_{\mp} & \mp \psi_{\pm} &      0 
             \end{array} 
    \right] , \\ 
    \hat{\nabla}_{T_2} 
   (\Theta_{\mp , 1} \ \Theta_{\pm , 2} \ \Theta_{\pm , 3} ) 
& =(\Theta_{\mp , 1} \ \Theta_{\pm , 2} \ \Theta_{\pm , 3} ) 
    \left[ \begin{array}{ccc} 
            0       & \pm  Z_{\pm}   & \mp  X_{\mp}   \\ 
        \pm Z_{\pm} &      0         & \mp \phi_{\mp} \\ 
        \mp X_{\mp} & \pm \phi_{\mp} &      0 
             \end{array} 
    \right] , 
\end{split} 
\label{hatnablans} 
\end{equation} 
where $W_{\pm}$, $X_{\pm}$, $Y_{\pm}$, $Z_{\pm}$ are as in \eqref{WXYZ} 
and $\phi_{\pm}$, $\psi_{\pm}$ are as in \eqref{phipsi}. 
Let $\hat{R}$ be the curvature tensor of $\hat{\nabla}$. 
Since $N$ is a space form of constant sectional curvature $L_0$, 
we obtain \eqref{hatR}, 
and computing the left side of \eqref{hatR} by \eqref{hatnablans}, 
we obtain  

\begin{pro}\label{pro:GCRns} 
The functions $W_{\pm}$, $X_{\pm}$, $Y_{\pm}$, $Z_{\pm}$ satisfy 
not only \eqref{W=XY=Z} but also 
\begin{equation} 
 W_{\mp} X_{\pm} +Y_{\pm} Z_{\mp} 
+L_0 e^{2\lambda} +(\phi_{\pm} )_u +(\psi_{\mp} )_v =0 
\label{gaussriccins} 
\end{equation} 
and \eqref{codazziWXYZ} with \eqref{phipsi}. 
\end{pro} 

\begin{rem} 
By direct computations as in \cite{AH}, 
we see that 
the system of two equations in \eqref{gaussriccins} is equivalent to 
the system of equations \eqref{gaussns}, \eqref{riccins}. 
\end{rem} 

Let $K$, $\nabla^{\perp}$, $R^{\perp}$ be as in 
Section~\ref{sect:Riem}. 
Then from \eqref{hatnablans} and \eqref{gaussriccins}, 
we obtain an analogue of Proposition~\ref{pro:WX+YZ=0} 
for space-like surfaces in $N$. 
In addition, noticing \eqref{gaussriccins}, 
we obtain analogues of 
Propositions~\ref{pro:nondegL=0}, \ref{pro:nondegLnot=0}, 
and Theorem~\ref{thm:csect}. 

\section{Time-like surfaces 
in neutral 4-dimensional space forms}\label{sect:nt} 

\setcounter{equation}{0} 

\subsection{The equations of Gauss, Codazzi and Ricci} 

Let $N$, $h$ be as in Section~\ref{sect:ns}. 
Let $M$ be a Lorentz surface 
and $F:M\longrightarrow N$ a time-like and conformal immersion 
of $M$ into $N$. 
Let $(u, v)$ be local coordinates of $M$ 
which are compatible with the paracomplex structure of $M$. 
Then the induced metric $g$ on $M$ by $F$ is represented 
as $g=e^{2\lambda} (du^2 -dv^2 )$ 
for a real-valued function $\lambda$. 
Let $\Tilde{\nabla}$, $T_1$, $T_2$ be as in Section~\ref{sect:ns}. 
Let $N_1$, $N_2$ be normal vector fields of $F$ satisfying 
\begin{equation} 
h(N_1 , N_1 )=-h(N_2 , N_2 )=e^{2\lambda} , \quad 
h(N_1 , N_2 )=0. 
\label{N1N2nt} 
\end{equation} 
Then we have \eqref{dt1t2} with 
\begin{equation} 
\begin{split} 
S & =\left[ \begin{array}{ccccc} 
         \lambda_u         & \lambda_v & -\alpha_1  &  \beta_1   & 1 \\ 
         \lambda_v         & \lambda_u &  \alpha_2  & -\beta_2   & 0 \\ 
         \alpha_1          & \alpha_2  &  \lambda_u &  \mu_1     & 0 \\ 
         \beta_1           & \beta_2   &  \mu_1     &  \lambda_u & 0 \\ 
         -L_0 e^{2\lambda} &  0        &   0        &   0        & 0 
              \end{array} 
     \right] , \\ 
T & =\left[ \begin{array}{ccccc} 
         \lambda_v &  \lambda_u        & -\alpha_2  &  \beta_2   & 0 \\ 
         \lambda_u &  \lambda_v        &  \alpha_3  & -\beta_3   & 1 \\ 
         \alpha_2  &  \alpha_3         &  \lambda_v &  \mu_2     & 0 \\ 
         \beta_2   &  \beta_3          &  \mu_2     &  \lambda_v & 0 \\ 
          0        &  L_0 e^{2\lambda} &   0        &   0        & 0 
              \end{array} 
    \right] , 
\end{split} 
\label{STnt} 
\end{equation} 
and $\alpha_k$, $\beta_k$ ($k=1, 2, 3$) and $\mu_l$ ($l=1, 2$) are 
real-valued functions. 
From \eqref{dt1t2}, 
we obtain $S_v -T_u =ST-TS$ with \eqref{STnt}. 
This is equivalent to 
the system of the equations of Gauss, Codazzi and Ricci. 
The equation of Gauss is given by 
\begin{equation} 
  \lambda_{uu} -\lambda_{vv} +L_0 e^{2\lambda} 
= \alpha_1 \alpha_3 -\beta_1 \beta_3 -\alpha^2_2 +\beta^2_2 . 
\label{gaussnt}
\end{equation}
The equations of Codazzi are given by 
\begin{equation} 
\begin{split}
  (\alpha_1 )_v -(\alpha_2 )_u 
& =\alpha_2 \lambda_u -\alpha_3 \lambda_v 
   +\beta_2 \mu_1     -\beta_1  \mu_2  , \\ 
   (\alpha_2 )_v -(\alpha_3 )_u 
& = \alpha_1 \lambda_u -\alpha_2 \lambda_v 
   +\beta_3  \mu_1     -\beta_2  \mu_2  , \\ 
   (\beta_1 )_v -(\beta_2 )_u 
& = \beta_2  \lambda_u -\beta_3  \lambda_v 
   +\alpha_2 \mu_1     -\alpha_1 \mu_2  , \\ 
   (\beta_2 )_v -(\beta_3 )_u 
& = \beta_1  \lambda_u -\beta_2  \lambda_v 
   +\alpha_3 \mu_1     -\alpha_2 \mu_2 .  
\end{split} 
\label{codazzint} 
\end{equation} 
The equation of Ricci is given by 
\begin{equation} 
  (\mu_1 )_v -(\mu_2 )_u = \alpha_1 \beta_2 -\alpha_2 \beta_1 
                          -\alpha_2 \beta_3 +\alpha_3 \beta_2 . 
\label{riccint} 
\end{equation}

Suppose that $N$ is oriented and 
that $(T_1 , N_1 , T_2 , N_2 )$ gives the orientation. 
We set 
\begin{equation*} 
e_1 :=\dfrac{1}{e^{\lambda}} T_1 , \quad 
e_2 :=\dfrac{1}{e^{\lambda}} N_1 , \quad 
e_3 :=\dfrac{1}{e^{\lambda}} T_2 , \quad 
e_4 :=\dfrac{1}{e^{\lambda}} N_2 . 
\end{equation*} 
Then $\Theta_{\mp , 1}$, $\Theta_{\pm , 2}$, $\Theta_{\pm , 3}$  
as in \eqref{Omega} form pseudo-orthonormal local frame fields 
of $\bigwedge^2_{\pm}\!F^*T\!N$ respectively. 
The \textit{time-like twistor lifts\/} of $F$ are 
sections of the time-like twistor spaces associated with $F^*T\!N$, 
and in this situation, 
these sections are locally given by $\Theta_{\pm , 2}$ respectively, 
and determined by $F$. 
Let $\nabla$, $\hat{\nabla}$ be as in Section~\ref{sect:Riem}. 
Then $\hat{\nabla}$ gives connections of $\bigwedge^2_{\pm}\!F^*T\!N$, 
and by \eqref{dt1t2} with \eqref{STnt}, we obtain 
\begin{equation} 
\begin{split} 
    \hat{\nabla}_{T_1}  
   (\Theta_{\mp , 1} \ \Theta_{\pm , 2} \ \Theta_{\pm , 3} ) 
& =(\Theta_{\mp , 1} \ \Theta_{\pm , 2} \ \Theta_{\pm , 3} ) 
    \left[ \begin{array}{ccc} 
            0         & W_{\pm} & \mp \psi_{\pm} \\ 
            W_{\pm}   & 0       &  -   Y_{\pm}   \\ 
       \mp \psi_{\pm} & Y_{\pm} &      0 
             \end{array} 
    \right] , \\ 
    \hat{\nabla}_{T_2} 
   (\Theta_{\mp , 1} \ \Theta_{\pm , 2} \ \Theta_{\pm , 3} ) 
& =(\Theta_{\mp , 1} \ \Theta_{\pm , 2} \ \Theta_{\pm , 3} ) 
    \left[ \begin{array}{ccc} 
            0         & \pm Z_{\pm} & \mp \phi_{\pm} \\ 
       \pm  Z_{\pm}   &     0       & \mp  X_{\pm} \\ 
       \mp \phi_{\pm} & \pm X_{\pm} &      0 
             \end{array} 
    \right] , 
\end{split} 
\label{hatnablant} 
\end{equation} 
where $W_{\pm}$, $X_{\pm}$, $Y_{\pm}$, $Z_{\pm}$ are as in \eqref{WXYZ} 
and $\phi_{\pm}$, $\psi_{\pm}$ are as in \eqref{phipsi}. 
In addition, 
noticing that $N$ is a space form of constant sectional curvature $L_0$, 
we obtain 

\begin{pro}\label{pro:GCRnt} 
The functions $W_{\pm}$, $X_{\pm}$, $Y_{\pm}$, $Z_{\pm}$ satisfy 
not only \eqref{W=XY=Z} but also 
\begin{equation} 
 W_{\pm} X_{\pm}  -Y_{\pm} Z_{\pm} 
+L_0 e^{2\lambda} +(\phi_{\pm} )_u -(\psi_{\pm} )_v =0 
\label{gaussriccint} 
\end{equation} 
and 
\begin{equation} 
\begin{split}
       (Y_{\pm} )_v \mp (X_{\pm} )_u 
& = \pm W_{\pm} \phi_{\pm} -Z_{\pm} \psi_{\pm} , \\ 
       (W_{\pm} )_v \mp (Z_{\pm} )_u 
& = \pm Y_{\pm} \phi_{\pm} -X_{\pm} \psi_{\pm} .  
\end{split} 
\label{codazziWXYZnt} 
\end{equation} 
\end{pro} 

\begin{rem} 
By direct computations as in \cite{AH}, 
we see that 
the system of two equations in \eqref{gaussriccint} is equivalent to 
the system of equations \eqref{gaussnt}, \eqref{riccint} and that 
the system of four equations in \eqref{codazziWXYZnt} is equivalent to 
the system of four equations in \eqref{codazzint}. 
\end{rem} 

Let $K$, $\nabla^{\perp}$, $R^{\perp}$ be as in 
Section~\ref{sect:Riem}. 
Then from \eqref{hatnablant} and \eqref{gaussriccint}, 
we obtain an analogue of Proposition~\ref{pro:WX+YZ=0} 
for time-like surfaces in $N$. 
In addition, noticing \eqref{gaussriccint} and \eqref{codazziWXYZnt}, 
we obtain analogues of 
Propositions~\ref{pro:nondegL=0}, \ref{pro:nondegLnot=0}. 

\subsection{Sections with fully light-like covariant derivative} 

Let $\Theta_{\varepsilon}$ be a section of 
a time-like twistor space $U_- \left( \bigwedge^2_{\varepsilon}\!F^*T\!N 
                               \right)$ 
in $\bigwedge^2_{\varepsilon}\!F^*T\!N$ ($\varepsilon \in \{ +, -\}$). 
Then $\Theta_{\varepsilon}$ is locally represented as 
$\Theta_{\varepsilon} 
=c_{\varepsilon , 1} \Theta_{-\varepsilon , 1} 
+c_{\varepsilon , 2} \Theta_{\varepsilon , 2} 
+c_{\varepsilon , 3} \Theta_{\varepsilon , 3}$, 
where $-+$, $--$ mean $-$, $+$ respectively, 
and $c_{\varepsilon , k}$ ($k=1, 2, 3$) are functions satisfying 
\begin{equation} 
 c^2_{\varepsilon , 1} 
-c^2_{\varepsilon , 2} 
-c^2_{\varepsilon , 3} 
=-1. 
\label{c} 
\end{equation} 
Since 
\begin{equation} 
 \hat{R} (T_1 , T_2 )\Theta_{\varepsilon} 
=\varepsilon L_0 e^{2\lambda} 
 (c_{\varepsilon , 3} \Theta_{-\varepsilon , 1} 
 +c_{\varepsilon , 1} \Theta_{\varepsilon , 3} ), 
\label{hatRTheta}
\end{equation} 
if $\Theta_{\varepsilon}$ is horizontal with respect to $\hat{\nabla}$ 
and if $L_0 \not= 0$, 
then $c_{\varepsilon , 1} =c_{\varepsilon , 3} =0$, 
and this means that $\Theta_{\varepsilon}$ is given by 
one of $\pm \Theta_{\varepsilon , 2}$. 
If $L_0 =0$, then $\hat{\nabla}$ is flat and therefore 
a local horizontal section can be constructed 
for an arbitrarily given initial value. 

Suppose that the covariant derivative $\hat{\nabla} \Theta_{\varepsilon}$ 
of $\Theta_{\varepsilon}$ is fully light-like, that is, 
$\hat{\nabla} \Theta_{\varepsilon}$ determines 
a light-like one-dimensional subspace of $\bigwedge^2_{\varepsilon}\!F^*T\!N$ 
at each point of $M$. 
Noticing $\hat{h} (\hat{\nabla} \Theta_{\varepsilon} , 
                                \Theta_{\varepsilon} )=0$ 
and that $\hat{\nabla} \Theta_{\varepsilon}$ is fully light-like, 
we obtain $\hat{h} (\hat{R} (T_1 , T_2 )\Theta_{\varepsilon} , 
                                        \Theta_{\varepsilon} )=0$. 
Since $\hat{h} (\hat{R} (T_1 , T_2 )\Theta_{\varepsilon} , 
                       \hat{\nabla} \Theta_{\varepsilon} )=0$, 
$\hat{R} (T_1 , T_2 )\Theta_{\varepsilon}$ 
and    $\hat{\nabla} \Theta_{\varepsilon}$ are linearly dependent. 
In particular, $\hat{R} (T_1 , T_2 )\Theta_{\varepsilon}$ is 
zero or light-like. 
Moreover, by \eqref{hatRTheta}, 
if $L_0 \not= 0$, 
then $c_{\varepsilon , 3} =\delta_{\varepsilon} c_{\varepsilon , 1}$ 
for $\delta_{\varepsilon} =1$ or $-1$. 
Therefore, in this case, from \eqref{c}, 
we obtain $c_{\varepsilon , 2} =1$ or $-1$. 

In the following, we suppose 
\begin{equation} 
\Theta_{\varepsilon} 
=c_{\varepsilon , 1} \Theta_{-\varepsilon , 1} 
+                    \Theta_{\varepsilon , 2} 
+\delta_{\varepsilon} 
 c_{\varepsilon , 1} \Theta_{\varepsilon , 3} . 
\label{Theta+} 
\end{equation} 
Then by the first relation in \eqref{hatnablant}, we obtain 
\begin{equation} 
  \hat{\nabla}_{T_1} \Theta_{\varepsilon} 
=(\Theta_{-\varepsilon , 1} \ 
  \Theta_{\varepsilon , 2}  \ 
  \Theta_{\varepsilon , 3} ) 
  \left[ \begin{array}{c} 
          (c_{\varepsilon , 1} )_u 
          -   \varepsilon   \delta_{\varepsilon} \psi_{\varepsilon} 
           c_{\varepsilon , 1}  
          +W_{\varepsilon} \\ 
          (W_{\varepsilon} -\delta_{\varepsilon}    Y_{\varepsilon} ) 
           c_{\varepsilon , 1} \\    
                            \delta_{\varepsilon} 
          (c_{\varepsilon , 1} )_u 
             -   \varepsilon \psi_{\varepsilon} 
           c_{\varepsilon , 1} 
          +Y_{\varepsilon} 
           \end{array} 
  \right] . 
\label{duTheta} 
\end{equation} 
Since $\hat{\nabla} \Theta_{\varepsilon}$ is zero or light-like, 
we have 
\begin{equation*} 
 ((c_{\varepsilon , 1} )_u 
  -   \varepsilon   \delta_{\varepsilon} \psi_{\varepsilon} 
   c_{\varepsilon , 1}  
  +W_{\varepsilon} )^2 
-( W_{\varepsilon} -\delta_{\varepsilon}    Y_{\varepsilon} )^2 
 c^2_{\varepsilon , 1} 
-(\delta_{\varepsilon} (c_{\varepsilon , 1} )_u 
  -   \varepsilon \psi_{\varepsilon} 
   c_{\varepsilon , 1} 
  +Y_{\varepsilon} )^2 
=0. 
\end{equation*} 
This yields 
\begin{equation} 
2(W_{\varepsilon} -\delta_{\varepsilon} Y_{\varepsilon} ) 
 (c_{\varepsilon , 1} )_u 
=(W_{\varepsilon} -\delta_{\varepsilon} Y_{\varepsilon} )^2 
  c^2_{\varepsilon , 1} 
+2\varepsilon \psi_{\varepsilon} 
 (\delta_{\varepsilon} W_{\varepsilon} -Y_{\varepsilon} ) 
  c_{\varepsilon , 1} 
- W^2_{\varepsilon} +Y^2_{\varepsilon} . 
\label{cu} 
\end{equation} 
Suppose $W_{\varepsilon} \not= \delta_{\varepsilon} Y_{\varepsilon}$. 
Then from \eqref{cu}, we have 
\begin{equation} 
2(c_{\varepsilon , 1} )_u 
=(W_{\varepsilon} -\delta_{\varepsilon} Y_{\varepsilon} ) 
  c^2_{\varepsilon , 1} 
+2\varepsilon \delta_{\varepsilon} \psi_{\varepsilon} 
  c_{\varepsilon , 1} 
-(W_{\varepsilon} +\delta_{\varepsilon} Y_{\varepsilon} ). 
\label{cu2} 
\end{equation} 
By the second relation in \eqref{hatnablant}, we obtain 
\begin{equation} 
  \hat{\nabla}_{T_2} \Theta_{\varepsilon} 
=(\Theta_{-\varepsilon , 1} \ 
  \Theta_{\varepsilon , 2}  \ 
  \Theta_{\varepsilon , 3} ) 
  \left[ \begin{array}{c} 
          (c_{\varepsilon , 1} )_v 
          -   \varepsilon   \delta_{\varepsilon} \phi_{\varepsilon} 
           c_{\varepsilon , 1}  
          +   \varepsilon Z_{\varepsilon} \\ 
              \varepsilon 
          (Z_{\varepsilon} -\delta_{\varepsilon}    X_{\varepsilon} ) 
           c_{\varepsilon , 1} \\    
                            \delta_{\varepsilon} 
          (c_{\varepsilon , 1} )_v 
             -\varepsilon \phi_{\varepsilon} 
           c_{\varepsilon , 1} 
          +   \varepsilon X_{\varepsilon} 
           \end{array} 
  \right] . 
\label{dvTheta} 
\end{equation} 
Since $\hat{\nabla} \Theta_{\varepsilon}$ is zero or light-like, 
we have 
\begin{equation} 
2\varepsilon (Z_{\varepsilon} -\delta_{\varepsilon} X_{\varepsilon} ) 
 (c_{\varepsilon , 1} )_v 
=(Z_{\varepsilon} -\delta_{\varepsilon} X_{\varepsilon} )^2 
  c^2_{\varepsilon , 1} 
+2                 \delta_{\varepsilon} \phi_{\varepsilon} 
 (Z_{\varepsilon} -\delta_{\varepsilon} X_{\varepsilon} ) 
  c_{\varepsilon , 1} 
- Z^2_{\varepsilon} +X^2_{\varepsilon} . 
\label{cv} 
\end{equation} 
Suppose $Z_{\varepsilon} \not= \delta_{\varepsilon} X_{\varepsilon}$. 
Then from \eqref{cv}, we have 
\begin{equation} 
2\varepsilon (c_{\varepsilon , 1} )_v 
=(Z_{\varepsilon} -\delta_{\varepsilon} X_{\varepsilon} ) 
  c^2_{\varepsilon , 1} 
+2                 \delta_{\varepsilon} \phi_{\varepsilon} 
  c_{\varepsilon , 1} 
-(Z_{\varepsilon} +\delta_{\varepsilon} X_{\varepsilon} ). 
\label{cv2} 
\end{equation} 
Applying \eqref{gaussriccint}, \eqref{codazziWXYZnt}, \eqref{cu2} and 
\eqref{cv2} to $(c_{\varepsilon , 1} )_{uv} =(c_{\varepsilon , 1} )_{vu}$, 
we obtain $L_0 e^{2\lambda} c_{\varepsilon , 1} =0$.  
Therefore, if $L_0 \not= 0$, then $c_{\varepsilon , 1} =0$, 
and then from \eqref{cu2} and \eqref{cv2}, 
we obtain $W_{\varepsilon} =-\delta_{\varepsilon} Y_{\varepsilon}$ 
and       $Z_{\varepsilon} =-\delta_{\varepsilon} X_{\varepsilon}$. 
Hence we obtain 

\begin{pro}\label{pro:Lnot=0} 
Suppose $W_{\varepsilon} \not= \delta_{\varepsilon} Y_{\varepsilon}$ 
and     $Z_{\varepsilon} \not= \delta_{\varepsilon} X_{\varepsilon}$. 
If $L_0 \not= 0$, 
then a local section $\Theta_{\varepsilon}$ 
of $U_- \left( \bigwedge^2_{\varepsilon}\!F^*T\!N \right)$ such that 
$\hat{\nabla} \Theta_{\varepsilon}$ is fully light-like is given by 
one of $\pm \Theta_{\varepsilon , 2}$ 
and then $W_{\varepsilon} =-\delta_{\varepsilon} Y_{\varepsilon}$ 
and      $Z_{\varepsilon} =-\delta_{\varepsilon} X_{\varepsilon}$. 
\end{pro} 

\begin{rem} 
Suppose $L_0 =0$. 
Then $\hat{\nabla}$ is flat and 
we can construct local sections 
of $U_- \left( \bigwedge^2_{\varepsilon}\!F^*T\!N \right)$ 
with fully light-like covariant derivative (see \cite{ando8}). 
\end{rem} 

By \eqref{duTheta} and \eqref{dvTheta}, we obtain 

\begin{pro}\label{pro:=0not=0} 
Suppose $L_0 \not= 0$. 
\begin{itemize} 
\item[{\rm (a)}]{If 
     $W_{\varepsilon}     = \delta_{\varepsilon} Y_{\varepsilon}$ and 
     $Z_{\varepsilon} \not= \delta_{\varepsilon} X_{\varepsilon}$ hold, 
then $c_{\varepsilon , 1} =  0$, 
     $W_{\varepsilon} =0$ 
and  $Z_{\varepsilon}     =-\delta_{\varepsilon} X_{\varepsilon} ;$} 
\item[{\rm (b)}]{if 
     $W_{\varepsilon} \not= \delta_{\varepsilon} Y_{\varepsilon}$ and 
     $Z_{\varepsilon}     = \delta_{\varepsilon} X_{\varepsilon}$ hold, 
then $c_{\varepsilon , 1} =  0$, 
     $Z_{\varepsilon} =0$ 
and  $W_{\varepsilon}     =-\delta_{\varepsilon} Y_{\varepsilon}$.} 
\end{itemize} 
\end{pro} 

Noticing Propositions~\ref{pro:Lnot=0}, \ref{pro:=0not=0},  
in order to study a section $\Theta_{\varepsilon}$ 
with fully light-like covariant derivative in the case of $L_0 \not= 0$, 
we can suppose 
\begin{equation} 
(Z_{\varepsilon} , W_{\varepsilon} )
   =\delta_{\varepsilon} 
(X_{\varepsilon} , Y_{\varepsilon} ) 
    \quad 
   (\varepsilon =+, -). 
\label{ZWXY} 
\end{equation} 
Then for any $\varepsilon \in \{ +, -\}$ 
and a function $c_{\varepsilon , 1}$ satisfying at least one of 
\begin{equation} 
           \delta_{\varepsilon} 
               (c_{\varepsilon , 1} )_u 
-\varepsilon \psi_{\varepsilon} 
                c_{\varepsilon , 1} 
               +Y_{\varepsilon} \not= 0, \quad 
           \delta_{\varepsilon} 
               (c_{\varepsilon , 1} )_v 
-\varepsilon \phi_{\varepsilon} 
                c_{\varepsilon , 1} 
+\varepsilon    X_{\varepsilon} \not= 0, 
\label{flc} 
\end{equation} 
the covariant derivative of $\Theta_{\varepsilon}$ as in \eqref{Theta+} is 
fully light-like. 
Applying \eqref{ZWXY} to \eqref{gaussriccint}, we obtain 
\begin{equation} 
L_0 e^{2\lambda} +(\phi_{\varepsilon} )_u -(\psi_{\varepsilon} )_v =0 
\quad (\varepsilon =+, -). 
\label{gaussriccint2} 
\end{equation} 
Applying \eqref{phipsi} to \eqref{gaussriccint2}, 
we see that $F$ satisfies $K\equiv L_0$ and $R^{\perp} \equiv 0$. 
The former condition is equivalent to 
\begin{equation} 
\lambda_{uu} -\lambda_{vv} +L_0 e^{2\lambda} =0 
\label{K=L0tn} 
\end{equation} 
and the latter condition is equivalent to the existence of 
a function $\gamma$ satisfying \eqref{gamma}. 
From \eqref{codazziWXYZnt}, 
we obtain 
\begin{equation} 
(Y_{\varepsilon} )_v -\varepsilon (X_{\varepsilon} )_u 
=\delta_{\varepsilon} 
(\phi_{\varepsilon} Y_{\varepsilon} -\psi_{\varepsilon} X_{\varepsilon} ) 
\quad (\varepsilon =+, -). 
\label{codazziWXYZnt2} 
\end{equation} 
From \eqref{W=XY=Z}, we have 
\begin{equation} 
\delta_+ Y_+ +\delta_- Y_- =X_+ +X_- , \quad 
\delta_+ X_+ +\delta_- X_- =Y_+ +Y_- . 
\label{deltaXY} 
\end{equation} 
From \eqref{deltaXY}, we have 
\begin{equation} 
(\delta_+ \delta_- -1)Y_+ =(\delta_- -\delta_+ )X_+ , \quad 
(\delta_+ \delta_- -1)Y_- =(\delta_+ -\delta_- )X_- . 
\label{deltaXY2} 
\end{equation} 

Suppose $\delta_- =-\delta_+$. 
Then from \eqref{ZWXY} and \eqref{deltaXY2}, 
we obtain $W_{\pm} =X_{\pm}$ and $Y_{\pm} =Z_{\pm}$. 
Therefore by the definitions of $W_{\pm}$, $X_{\pm}$, $Y_{\pm}$, $Z_{\pm}$ 
in \eqref{WXYZ}, $F$ has zero mean curvature vector. 

Suppose $\delta_- =\delta_+$ ($=:\delta$). 
Then from \eqref{phipsi} and \eqref{codazziWXYZnt2}, we have 
\begin{equation} 
\begin{split} 
(Y_+ )_v -(X_+ )_u & =\delta ((\lambda_u -\gamma_v )Y_+ 
                             -(\lambda_v -\gamma_u )X_+ ), \\  
(Y_- )_v +(X_- )_u & =\delta ((\lambda_u +\gamma_v )Y_- 
                             -(\lambda_v +\gamma_u )X_- ). 
\end{split} 
\label{codazziWXYZnt3} 
\end{equation} 
From \eqref{deltaXY}, we have 
\begin{equation} 
Y_+ +Y_- =\delta (X_+ +X_- ). 
\label{Y+-X+-} 
\end{equation} 
In this situation, 
the mean curvature vector of $F$ does not necessarily vanish, 
and $F$ has zero mean curvature vector 
if and only if $Y_{\pm} =\delta X_{\pm}$. 

Let $\lambda$ be a function of two variables $u$, $v$ 
satisfying \eqref{K=L0tn}. 
Let $\gamma$ be a function of $u$, $v$. 
Let $X_{\pm}$, $Y_{\pm}$ be functions of $u$, $v$ 
satisfying \eqref{codazziWXYZnt3} and \eqref{Y+-X+-} 
for $\delta \in \{ 1, -1\}$. 
Let $W_{\pm}$, $Z_{\pm}$ be defined as in \eqref{ZWXY} 
with $\delta_+ =\delta_- =\delta$. 
Let $\phi_{\pm}$, $\psi_{\pm}$ be defined as in \eqref{phipsi} 
with \eqref{gamma}. 
Then these functions satisfy 
\eqref{W=XY=Z}, \eqref{gaussriccint} and \eqref{codazziWXYZnt}. 
Therefore these functions give an immersion into $N$ 
such that for any $\varepsilon \in \{ +, -\}$ 
and any function $c_{\varepsilon , 1}$ satisfying 
at least one of the two relations in \eqref{flc}, 
the covariant derivative of $\Theta_{\varepsilon}$ as in \eqref{Theta+} is 
fully light-like. 
This immersion is unique up to an isometry of $N$. 

Hence we obtain 

\begin{thm}\label{thm:tnfl} 
Suppose $L_0 \not= 0$. 
Then for each $\varepsilon \in \{ +, -\}$, 
the covariant derivative of $\Theta_{\varepsilon}$ is fully light-like 
if and only if $F$ satisfies \eqref{ZWXY} 
and $c_{\varepsilon , 1}$ satisfies 
at least one of the two relations in \eqref{flc}. 
If $F$ satisfies \eqref{ZWXY}, 
then 
\begin{itemize} 
\item[{\rm (a)}]{$K\equiv L_0$ and $R^{\perp} \equiv 0$,} 
\item[{\rm (b)}]{one of the following holds\/$:$ 
\begin{itemize} 
\item[{\rm (b1)}]{$F$ has zero mean curvature vector\/$;$} 
\item[{\rm (b2)}]{$\delta_+ =\delta_- (=\delta )$, and 
$X_{\pm}$, $Y_{\pm}$ satisfy \eqref{codazziWXYZnt3} and \eqref{Y+-X+-}.}  
\end{itemize}} 
\end{itemize}  
In addition, 
for a function $\lambda$ of two variables $u$, $v$ 
satisfying \eqref{K=L0tn}, 
a function $\gamma$ of $u$, $v$ and 
functions $X_{\pm}$, $Y_{\pm}$ of $u$, $v$ 
satisfying \eqref{codazziWXYZnt3} and \eqref{Y+-X+-} 
with $\delta \in \{ 1, -1\}$, 
an immersion into $N$ with $L_0 \not= 0$, 
\eqref{phipsi} and \eqref{ZWXY} is constructed 
and unique up to an isometry of $N$. 
\end{thm} 

\begin{rem}\label{rem:ttlflcd}  
The time-like twistor lifts of $F$ have 
fully light-like covariant derivative 
if and only if $F$ satisfies \eqref{ZWXY} 
so that the vectors in \eqref{ZWXY} are nonzero. 
\end{rem} 

\begin{rem}\label{rem:ldc}  
Suppose that the immersion $F$ satisfies \eqref{ZWXY}. 
Then the second fundamental form $\sigma$ of $F$ satisfies 
the linearly dependent condition, and 
$F$ satisfies (b2) in Theorem~\ref{thm:tnfl} if and only if 
$\sigma$ satisfies the light-like linearly dependent condition, i.e., 
\eqref{ldc} with $\theta =\pi /4$ or $-\pi /4$. 
This condition just means that 
a light-like normal vector field of $F$ is parallel (\cite{AH}). 
\end{rem} 

\begin{rem} 
By \eqref{WXYZ}, 
the mean curvature vector of $F$ vanishes 
if and only if $F$ satisfies 
\begin{equation} 
(W_{\pm} , Z_{\pm} )=(X_{\pm} , Y_{\pm} ). 
\label{WZXYtn} 
\end{equation} 
Therefore $F$ satisfies \eqref{ZWXY} and \eqref{WZXYtn} 
if and only if 
\begin{equation} 
W_{\pm} =X_{\pm} =\delta_{\pm} Y_{\pm} =\delta_{\pm} Z_{\pm} , 
\label{H=0} 
\end{equation} 
and these functions are not necessarily zero. 
Refer to \cite{ando5}, \cite{ando9}, \cite{AH} 
for characterizations of time-like surfaces in $N$ 
with zero mean curvature vector 
such that their time-like twistor lifts have 
zero or light-like covariant derivative. 
\end{rem} 

\begin{rem}\label{rem:pcqd} 
Let $F:M\longrightarrow N$ be a time-like and conformal immersion. 
Let $\check{w} :=u+jv$ be a local paracomplex coordinate of $M$ 
and set $\check{\partial} 
       :=\partial /\partial \check{w} 
  =(1/2)(\partial /\partial u +j\partial /\partial v)$, 
where $j$ is the paraimaginary unit of paracomplex numbers. 
Let $\check{Q}$ be a paracomplex quartic differential on $M$ defined by 
\begin{equation} 
\check{Q} =h(\sigma (\check{\partial} , \check{\partial} ), 
             \sigma (\check{\partial} , \check{\partial} ))d\check{w}^4 
\label{QNt} 
\end{equation} 
for the second fundamental form $\sigma$ of $F$. 
Then we obtain 
\begin{equation} 
h(\sigma (\check{\partial} , \check{\partial} ), 
  \sigma (\check{\partial} , \check{\partial} ))
= \dfrac{e^{2\lambda}}{16} 
   (          \check{p}  +2j          \check{q}  ) 
   (\overline{\check{p}} +2j\overline{\check{q}} ), 
\label{QNtc} 
\end{equation} 
where 
\begin{equation}
\check{p} :=\alpha_1 +\alpha_3 +j(\beta_1 +\beta_3 ), \quad 
\check{q} :=\alpha_2           +j \beta_2 . 
\label{checkpq} 
\end{equation} 
Therefore $\check{Q}$ is zero or null if and only if 
one of $\check{p} \pm 2\check{q}$ is zero or null. 
Noticing 
\begin{equation} 
\check{p} +2\delta \check{q} 
=Y_+ -Z_- +\delta (W_+ +W_- ) +j(W_+ -X_- +\delta (Y_+ +Y_- )) 
\label{ppm2qNt} 
\end{equation} 
for $\delta \in \{ 1, -1\}$, 
we see that if $F$ satisfies \eqref{ZWXY}, 
then $\check{Q}$ is zero or null. 
If the second fundamental form $\sigma$ of $F$ satisfies 
the light-like linearly dependent condition, 
then by \eqref{QNtc} and \eqref{checkpq}, $\check{Q}$ is zero. 
\end{rem} 

\subsection{Time-like twistor lifts which are horizontal 
along integral curves of a light-like one-dimensional distribution} 

In the following, 
we study a time-like surface in $N$ such that 
for any $\varepsilon \in \{ +, -\}$, 
the covariant derivative of $\Theta_{\varepsilon , 2}$ 
along integral curves of a light-like one-dimensional distribution 
on the surface vanishes. 
By \eqref{hatnablant}, we have 
\begin{equation*} 
\hat{\nabla}_{T_1 +T_2} \Theta_{\varepsilon , 2} 
=(W_{\varepsilon} +\varepsilon Z_{\varepsilon} )\Theta_{-\varepsilon , 1} 
+(Y_{\varepsilon} +\varepsilon X_{\varepsilon} )\Theta_{\varepsilon , 3} . 
\label{nabladu+dvT2} 
\end{equation*} 
Therefore $\hat{\nabla}_{T_1 +T_2} \Theta_{\varepsilon , 2}$ vanishes 
if and only if $F$ satisfies 
\begin{equation} 
              (W_{\varepsilon} , X_{\varepsilon} ) 
=-\varepsilon (Z_{\varepsilon} , Y_{\varepsilon} ) \quad 
(\varepsilon =+, -). 
\label{WXZY} 
\end{equation} 
Applying \eqref{WXZY} to \eqref{gaussriccint}, 
we obtain \eqref{gaussriccint2}, 
and therefore we obtain $K\equiv L_0$ and $R^{\perp} \equiv 0$. 
Applying \eqref{WXZY} to \eqref{W=XY=Z}, 
we obtain 
\begin{equation} 
W_+ =X_+ =-Y_+ =-Z_+ , \quad 
W_- =X_- = Y_- = Z_- , 
\label{W+W-} 
\end{equation} 
which is \eqref{H=0} with $(\delta_+ , \delta_- )=(-1, 1)$. 
Applying \eqref{phipsi} with \eqref{gamma} and \eqref{W+W-} 
to \eqref{codazziWXYZnt}, 
we obtain 
\begin{equation*} 
\begin{split} 
        (\partial_u +\partial_v )X_+ 
& =-X_+ (\partial_u +\partial_v )(\lambda -\gamma ), \\ 
        (\partial_u +\partial_v )X_- 
& =-X_- (\partial_u +\partial_v )(\lambda +\gamma ), 
\end{split} 
\label{dudvlg} 
\end{equation*} 
which is equivalent to 
\begin{equation} 
(\partial_u +\partial_v )(X_+ e^{\lambda -\gamma} )=0, \quad 
(\partial_u +\partial_v )(X_- e^{\lambda +\gamma} )=0. 
\label{dudvlg2} 
\end{equation} 
If $X_{\pm}$, $\lambda$, $\gamma$ satisfy \eqref{dudvlg2}, 
then there exist functions $\xi_{\pm}$ of one variable 
satisfying 
\begin{equation} 
   X_{\varepsilon} 
=\xi_{\varepsilon} (u-v)e^{-\lambda +\varepsilon \gamma} \quad 
     (\varepsilon =+, -). 
\label{xipm} 
\end{equation} 

Let $\lambda$ be a function of two variables $u$, $v$ 
satisfying \eqref{K=L0tn}. 
Let $\gamma$ be a function of $u$, $v$. 
Let $X_{\pm}$ be functions of $u$, $v$ 
satisfying \eqref{xipm} for functions $\xi_{\pm}$ of one variable. 
Let $W_{\pm}$, $Y_{\pm}$, $Z_{\pm}$ be defined as in \eqref{W+W-}. 
Let $\phi_{\pm}$, $\psi_{\pm}$ be defined as in \eqref{phipsi} 
with \eqref{gamma}. 
Then these functions satisfy 
\eqref{W=XY=Z}, \eqref{gaussriccint} and \eqref{codazziWXYZnt}. 
Therefore these functions give an immersion into $N$ 
satisfying $\hat{\nabla}_{T_1 +T_2} \Theta_{\varepsilon , 2} =0$, 
which is unique up to an isometry of $N$. 

We can similarly study 
the condition $\hat{\nabla}_{T_1 -T_2} \Theta_{\varepsilon , 2} =0$. 
This condition is equivalent to 
\begin{equation} 
             (W_{\varepsilon} , X_{\varepsilon} ) 
=\varepsilon (Z_{\varepsilon} , Y_{\varepsilon} ) \quad 
(\varepsilon =+, -), 
\label{WXZY-} 
\end{equation} 
instead of \eqref{WXZY}. 
Therefore we have 
\begin{equation} 
W_+ =X_+ = Y_+ = Z_+ , \quad 
W_- =X_- =-Y_- =-Z_- , 
\label{W+W--} 
\end{equation} 
instead of \eqref{W+W-}, 
and \eqref{W+W--} is \eqref{H=0} with $(\delta_+ , \delta_- )=(1, -1)$. 
In addition, we have 
\begin{equation} 
   X_{\varepsilon} 
=\xi_{\varepsilon} (u+v)e^{-\lambda -\varepsilon \gamma} \quad 
     (\varepsilon =+, -) 
\label{xipm-} 
\end{equation} 
for functions $\xi_{\pm}$ of one variable, instead of \eqref{xipm}. 

A time-like surface in $N$ satisfies \eqref{W+W-} or \eqref{W+W--} 
if and only if the surface has zero mean curvature vector 
so that the shape operator of any normal vector field is 
zero or light-like. 

Hence we obtain 

\begin{thm}\label{thm:ldcdz} 
For a time-like and conformal immersion $F:M\longrightarrow N$, 
the following are mutually equivalent\/$:$ 
\begin{itemize} 
\item[{\rm (A)}]{there exists 
a light-like one-dimensional distribution $\mathcal{D}$ on $M$ 
such that for any $\varepsilon \in \{ +, -\}$, 
the covariant derivatives of $\Theta_{\varepsilon , 2}$ 
along integral curves of $\mathcal{D}$ vanish\/$;$} 
\item[{\rm (B)}]{$F$ has zero mean curvature vector so that 
the shape operator of any normal vector field is 
zero or light-like\/$;$} 
\item[{\rm (C)}]{$F$ satisfies \eqref{W+W-} or \eqref{W+W--}.} 
\end{itemize} 
If $F$ satisfies \eqref{W+W-} or \eqref{W+W--}, then 
\begin{itemize} 
\item[{\rm (a)}]{$K\equiv L_0$ and $R^{\perp} \equiv 0$,} 
\item[{\rm (b)}]{$X_{\pm}$ are represented as in 
either \eqref{xipm} or \eqref{xipm-}.} 
\end{itemize} 
In addition, 
for a function $\lambda$ satisfying \eqref{K=L0tn}, 
a function $\gamma$ and 
functions $X_{\pm}$ as in \eqref{xipm} or \eqref{xipm-} 
with functions $\xi_{\pm}$ of one variable, 
an immersion into $N$ with the above {\rm (A)} is constructed 
by \eqref{W+W-} or \eqref{W+W--}, 
and unique up to an isometry of $N$. 
\end{thm} 

\begin{rem} 
We can refer to \cite{ando5} 
for a characterization of time-like surfaces in $N$ 
as in (B) of Theorem~\ref{thm:ldcdz}. 
\end{rem} 

\subsection{Horizontality of the paracomplex lift} 

Let $\check{\mbox{\boldmath{$C$}}}$ be 
the set of paracomplex numbers. 
Let $\check{\Theta}$ be a section of 
the paracomplexification $\bigwedge^2\!F^*T\!N \otimes 
                          \check{\mbox{\boldmath{$C$}}}$ 
of $\bigwedge^2\!F^*T\!N$ which is locally given by 
\begin{equation} 
\begin{split} 
     \check{\Theta} 
:= & \dfrac{1}{\sqrt{2}} (e_1 \wedge e_3 +je_4 \wedge e_2 ) \\ 
 = & \dfrac{1}{2} ( \Theta_{+, 2} +\Theta_{-, 2} 
     +j            (\Theta_{+, 2} -\Theta_{-, 2} )),  
\end{split} 
\label{Thetant} 
\end{equation} 
where $j$ denotes the paraimaginary unit of paracomplex numbers 
as in Remark~\ref{rem:pcqd}. 
Then $\check{\Theta}$ in \eqref{Thetant} is determined by the immersion $F$ 
and the orientations of $M$, $N$. 
We call $\check{\Theta}$ the \textit{paracomplex lift\/}. 
Let $\check{w} :=u+jv$ be a local paracomplex coordinate of $M$. 
Then we can naturally consider 
the covariant derivative $\hat{\nabla}_{\check{\partial}} \check{\Theta}$ 
of $\check{\Theta}$ by $\check{\partial}$ and we have 
\begin{equation} 
\begin{split} 
    \hat{\nabla}_{\check{\partial}} \check{\Theta} 
= & \dfrac{1}{4} (  \hat{\nabla}_{T_1} (\Theta_{+, 2} +\Theta_{-, 2} ) 
                   +\hat{\nabla}_{T_2} (\Theta_{+, 2} -\Theta_{-, 2} ) \\ 
  & \quad          + j 
                  ( \hat{\nabla}_{T_1} (\Theta_{+, 2} -\Theta_{-, 2} ) 
                   +\hat{\nabla}_{T_2} (\Theta_{+, 2} +\Theta_{-, 2} ))).  
\end{split} 
\label{hndTcheck} 
\end{equation} 
Whether $\hat{\nabla}_{\check{\partial}} \check{\Theta}$ vanishes 
does not depend on the choice of $\check{w}$, and 
from \eqref{hndTcheck}, 
we see that $\hat{\nabla}_{\check{\partial}} \check{\Theta}$ vanishes 
if and only if the following hold: 
\begin{equation} 
  \hat{\nabla}_{T_1} \Theta_{+, 2} 
=-\hat{\nabla}_{T_2} \Theta_{+, 2} , \quad 
  \hat{\nabla}_{T_1} \Theta_{-, 2} 
= \hat{\nabla}_{T_2} \Theta_{-, 2} . 
\label{hndT0check} 
\end{equation} 
By \eqref{hatnablant}, \eqref{hndT0check} is equivalent to 
\begin{equation} 
(W_{\pm} , X_{\pm} )=-(Z_{\pm} , Y_{\pm} ). 
\label{WX-ZYcheck} 
\end{equation} 
Suppose \eqref{WX-ZYcheck}. 
Then by \eqref{gaussricci}, 
we obtain $K\equiv L_0$ and $R^{\perp} \equiv 0$. 
Applying \eqref{WX-ZYcheck} to \eqref{codazziWXYZnt}, we obtain 
\begin{equation} 
\begin{split} 
            (\partial_u \pm \partial_v )X_{\pm} 
& =-W_{\pm} (\partial_u \pm \partial_v )(\lambda -\gamma ), \\ 
            (\partial_u \pm \partial_v )W_{\pm} 
& =-X_{\pm} (\partial_u \pm \partial_v )(\lambda -\gamma ), 
\end{split} 
\label{dupmdvcheck} 
\end{equation} 
which is equivalent to 
\begin{equation*} 
\begin{split} 
(\partial_u +\partial_v )((W_+ \pm X_+ )e^{\pm (\lambda -\gamma )} ) & =0, \\ 
(\partial_u -\partial_v )((W_- \pm X_- )e^{\pm (\lambda -\gamma )} ) & =0. 
\end{split} 
\label{dupmdvcheck2} 
\end{equation*} 
Therefore there exist 
functions $\xi_{\pm}  =       \xi_{\pm}  (u\mp v)$, 
   $\Tilde{\xi}_{\pm} =\Tilde{\xi}_{\pm} (u\mp v)$ satisfying 
\begin{equation} 
       \xi_{\pm}  (u\mp v)=(W_{\pm} +X_{\pm} )e^{\lambda  -\gamma} , \quad 
\Tilde{\xi}_{\pm} (u\mp v)=(W_{\pm} -X_{\pm} )e^{-\lambda +\gamma} . 
\label{xiTildexi} 
\end{equation} 
From \eqref{xiTildexi}, we obtain 
\begin{equation*} 
\begin{split} 
W_{\pm} & =\dfrac{1}{2} (       \xi_{\pm}  (u\mp v)e^{-\lambda +\gamma} 
                        +\Tilde{\xi}_{\pm} (u\mp v)e^{\lambda  -\gamma} ), \\ 
X_{\pm} & =\dfrac{1}{2} (       \xi_{\pm}  (u\mp v)e^{-\lambda +\gamma} 
                        -\Tilde{\xi}_{\pm} (u\mp v)e^{\lambda  -\gamma} ). 
\end{split} 
\end{equation*} 
Then by \eqref{W=XY=Z}, we see that $\Tilde{\xi}_{\pm}$ are constant. 
Therefore there exist real numbers $r_{\pm}$ satisfying 
\begin{equation} 
\begin{split} 
W_{\pm} & =\dfrac{1}{2} ( \xi_{\pm} (u\mp v)e^{-\lambda +\gamma} 
                         +  r_{\pm}         e^{\lambda  -\gamma} ), \\ 
X_{\pm} & =\dfrac{1}{2} ( \xi_{\pm} (u\mp v)e^{-\lambda +\gamma} 
                         -  r_{\pm}         e^{\lambda  -\gamma} ). 
\end{split} 
\label{WpmXpmcheck} 
\end{equation} 

Let $\lambda$ be a function of two variables $u$, $v$ 
satisfying \eqref{K=L0tn}. 
Let $\gamma$ be a function of $u$, $v$. 
Let $W_{\pm}$, $X_{\pm}$ be functions of $u$, $v$ 
satisfying \eqref{WpmXpmcheck} 
for functions $\xi_{\pm}$ of one variable and real numbers $r_{\pm}$. 
Let $Y_{\pm}$, $Z_{\pm}$ be defined as in \eqref{WX-ZYcheck}. 
Let $\phi_{\pm}$, $\psi_{\pm}$ be defined as in \eqref{phipsi} 
with \eqref{gamma}. 
Then these functions satisfy 
\eqref{W=XY=Z}, \eqref{gaussriccint} and \eqref{codazziWXYZnt}. 
Therefore these functions give an immersion into $N$ 
satisfying $\hat{\nabla}_{\check{\partial}} \check{\Theta} =0$, 
which is unique up to an isometry of $N$. 

We can similarly study 
the condition $\hat{\nabla}_{\overline{\check{\partial}}} \check{\Theta} =0$ 
for   $\overline{\check{\partial}} 
     :=\partial /\partial \overline{\check{w}} 
=(1/2)(\partial /\partial u -j\partial /\partial v)$. 
This condition is equivalent to 
\begin{equation} 
(W_{\pm} , X_{\pm} )=(Z_{\pm} , Y_{\pm} ), 
\label{WX-ZYcheck-} 
\end{equation} 
instead of \eqref{WX-ZYcheck}. 
Suppose $\hat{\nabla}_{\overline{\check{\partial}}} \check{\Theta} =0$. 
Then applying \eqref{WX-ZYcheck-} to \eqref{codazziWXYZnt}, we obtain 
\begin{equation} 
\begin{split} 
            (\partial_u \mp \partial_v )X_{\pm} 
& =-W_{\pm} (\partial_u \mp \partial_v )(\lambda +\gamma ), \\ 
            (\partial_u \mp \partial_v )W_{\pm} 
& =-X_{\pm} (\partial_u \mp \partial_v )(\lambda +\gamma ), 
\end{split} 
\label{dupmdvcheck-} 
\end{equation} 
instead of \eqref{dupmdvcheck}, and 
\begin{equation} 
\begin{split} 
W_{\pm} & =\dfrac{1}{2} ( \xi_{\pm} (u\mp v)e^{-\lambda -\gamma} 
                         +  r_{\pm}         e^{\lambda  +\gamma} ), \\ 
X_{\pm} & =\dfrac{1}{2} ( \xi_{\pm} (u\mp v)e^{-\lambda -\gamma} 
                         -  r_{\pm}         e^{\lambda  +\gamma} ), 
\end{split} 
\label{WpmXpmcheck-} 
\end{equation} 
instead of \eqref{WpmXpmcheck}. 

Hence we obtain 

\begin{thm}\label{thm:pcsect} 
The paracomplex lift $\check{\Theta}$ defined in \eqref{Thetant} satisfies 
either $\hat{\nabla}_{\check{\partial}} \check{\Theta} =0$ 
or     $\hat{\nabla}_{\overline{\check{\partial}}} \check{\Theta} =0$ 
if and only if the immersion $F$ satisfies 
\eqref{WX-ZYcheck} or \eqref{WX-ZYcheck-}. 
If $F$ satisfies \eqref{WX-ZYcheck} or \eqref{WX-ZYcheck-}, 
then 
\begin{itemize} 
\item[{\rm (a)}]{$K\equiv L_0$ and $R^{\perp} \equiv 0$,} 
\item[{\rm (b)}]{$W_{\pm}$, $X_{\pm}$ are represented as in 
either \eqref{WpmXpmcheck} or \eqref{WpmXpmcheck-} 
for functions $\xi_{\pm}$ of one variable and real numbers $r_{\pm}$.} 
\end{itemize} 
In addition, 
for a function $\lambda$ of two variables $u$, $v$ satisfying \eqref{K=L0tn}, 
a function $\gamma$ of $u$, $v$ and 
functions $W_{\pm}$, $X_{\pm}$ as in the above {\rm (b)}, 
an immersion into $N$ satisfying 
either $\hat{\nabla}_{\check{\partial}} \check{\Theta} =0$ 
or     $\hat{\nabla}_{\overline{\check{\partial}}} \check{\Theta} =0$ 
is constructed by \eqref{WX-ZYcheck} or \eqref{WX-ZYcheck-} 
and unique up to an isometry of $N$. 
\end{thm} 

\begin{rem}\label{rem:ldc2} 
Suppose that $F$ satisfies \eqref{WX-ZYcheck} or \eqref{WX-ZYcheck-}. 
Then the second fundamental form of $F$ satisfies 
the linearly dependent condition if and only if 
one of the following holds: 
\begin{itemize} 
\item[{\rm (i)}]{$\alpha_1 =\alpha_3$, 
which means that $F$ has zero mean curvature vector;} 
\item[{\rm (ii)}]{$\alpha_1 +\alpha_3 =0$ and $\alpha_2 =0$, 
which mean $\xi_{\pm} \equiv 0$.} 
\end{itemize} 
In any case, 
$\sigma$ satisfies the light-like linearly dependent condition, 
which is equivalent to the condition that 
a light-like normal vector field of $F$ is parallel (\cite{AH}). 
\end{rem} 

\begin{rem}\label{rem:+=-} 
Suppose that $F$ satisfies \eqref{WX-ZYcheck} or \eqref{WX-ZYcheck-}. 
Then the following are mutually equivalent: 
\begin{itemize} 
\item[{\rm (a)}]{$F$ has zero mean curvature vector;} 
\item[{\rm (b)}]{$F$ satisfies \eqref{H=0} with $\delta_+ =\delta_-$;} 
\item[{\rm (c)}]{$r_{\pm} =0$.} 
\end{itemize} 
In addition, if $F$ satisfies one of these conditions, 
then the second fundamental form of $F$ satisfies 
the light-like linearly dependent condition and therefore 
a light-like normal vector field of $F$ is parallel. 
\end{rem} 

\begin{rem} 
Suppose that $F$ has zero mean curvature vector. 
Then $F$ satisfies \eqref{WX-ZYcheck} or \eqref{WX-ZYcheck-} 
if and only if a light-like normal vector field of $F$ is parallel. 
\end{rem} 

\begin{rem} 
The conformal Gauss map of a time-like surface of Willmore type 
with zero paraholomorphic quartic differential 
in a $3$-dimensional Lorentzian space form 
or a $3$-dimensional space form with signature $(+, -, -)$ defines 
a time-like immersion into $S^4_2$ or $H^4_2$ 
with zero mean curvature vector and 
a parallel light-like normal vector field (\cite{ando4}). 
\end{rem} 

\begin{rem}\label{rem:pclQ}  
If $F$ satisfies \eqref{WX-ZYcheck} or \eqref{WX-ZYcheck-}, 
then noticing \eqref{ppm2qNt}, we see that  
the paracomplex quartic differential $\check{Q}$ on $M$ 
defined as in \eqref{QNt} is zero or null. 
\end{rem} 

\section{The complexification of the two-fold exterior power 
of an oriented Lorentzian vector bundle of rank 4}\label{sect:cplx} 

\setcounter{equation}{0} 

Let $X$ be an oriented real vector space of dimension $4$. 
Then the two-fold exterior power $\bigwedge^2\!X$ of $X$ is 
a $6$-dimensional vector space. 
Let $h$ be a Lorentzian inner product of $X$ 
and $*$ Hodge's $*$-operator 
with respect to $h$ and the orientation of $X$. 
Then $*$ gives a linear transformation of $\bigwedge^2\!X$ 
satisfying $*^2 =*\circ *=-{\rm Id}_{\bigwedge^2\!X}$ and 
\begin{equation} 
\begin{split} 
*(e_1 \wedge e_2 ) & =-e_3 \wedge e_4 , \\ 
*(e_1 \wedge e_3 ) & =-e_4 \wedge e_2 , \\ 
*(e_1 \wedge e_4 ) & = e_2 \wedge e_3 
\end{split} 
\label{star} 
\end{equation} 
for an ordered pseudo-orthonormal basis $(e_1 , e_2 , e_3 , e_4 )$ of $X$ 
giving the orientation of $X$ and satisfying 
\begin{equation} 
\begin{split} 
& h(e_k , e_k )= 1  \ (k=1, 2, 3), \quad 
  h(e_4 , e_4 )=-1, \\ 
& h(e_k , e_l )= 0  \ (k\not= l). 
\end{split} 
\label{pob}
\end{equation} 
Therefore 
the complexification $\bigwedge^2\!X\otimes \mbox{\boldmath{$C$}}$ 
of $\bigwedge^2\!X$ is decomposed into 
two subspaces $\bigwedge^2_{\pm}\!X$ of complex dimension $3$ 
corresponding to the eigenvalues $\pm \sqrt{-1}$ of $*$ respectively, 
and \eqref{star} implies that $\bigwedge^2_+\!X$ is spanned by 
\begin{equation} 
\begin{split} 
\Theta_1 := & 
\dfrac{1}{\sqrt{2}} (e_1 \wedge e_2 +\sqrt{-1}  e_3 \wedge e_4 ), \\ 
\Theta_2 := & 
\dfrac{1}{\sqrt{2}} (e_1 \wedge e_3 +\sqrt{-1}  e_4 \wedge e_2 ), \\ 
\Theta_3 := & 
\dfrac{1}{\sqrt{2}} ( \sqrt{-1} e_1 \wedge e_4 +e_2 \wedge e_3 ), 
\end{split} 
\label{Omega+X} 
\end{equation} 
and that $\bigwedge^2_-\!X$ is spanned by 
\begin{equation} 
\begin{split} 
\overline{\Theta_1} := & 
\dfrac{1}{\sqrt{2}} (e_1 \wedge e_2 -\sqrt{-1}  e_3 \wedge e_4 ), \\ 
\overline{\Theta_2} := & 
\dfrac{1}{\sqrt{2}} (e_1 \wedge e_3 -\sqrt{-1}  e_4 \wedge e_2 ), \\ 
\overline{\Theta_3} := & 
\dfrac{1}{\sqrt{2}} (-\sqrt{-1} e_1 \wedge e_4 +e_2 \wedge e_3 ). 
\end{split} 
\label{Omega-X} 
\end{equation} 

\begin{rem} 
For an element $P\in SO(3,1)$ 
and an ordered pseudo-orthonormal basis $e=(e_1 , e_2 , e_3 , e_4 )$ 
of $X$ with \eqref{pob} giving the orientation, 
let $T_P$ be a linear transformation of $X$ given by $T_P e=eP$. 
Then $T_P$ preserves $h$. 
We set $\Tilde{T}_P (e_i \wedge e_j ):=T_P e_i \wedge T_P e_j$. 
For an element $\Theta =\sum_{1\leq i<j\leq 4} c_{ij} e_i \wedge e_j$ 
of $\bigwedge^2\!X\otimes \mbox{\boldmath{$C$}}$, we set 
\begin{equation*} 
  \Tilde{T}_{P} (\Theta ) 
:=\sum_{1\leq i<j\leq 4} c_{ij} \Tilde{T}_P (e_i \wedge e_j ). 
\end{equation*} 
Then we obtain a linear transformation $\Tilde{T}_{P}$ 
of $\bigwedge^2\!X\otimes \mbox{\boldmath{$C$}}$. 
We set 
\begin{equation*} 
\begin{split} 
& P_{1, 1} :=\left[ 
             \begin{array}{cccc} 
              \cos \theta & -\sin \theta & 0 & 0 \\ 
              \sin \theta &  \cos \theta & 0 & 0 \\ 
               0          &   0          & 1 & 0 \\ 
               0          &   0          & 0 & 1 
               \end{array} 
             \right] , \ 
  P_{1, 2} :=\left[ 
             \begin{array}{cccc} 
              1          &  0          &   0           &   0 \\ 
              0          &  1          &   0           &   0 \\ 
              0          &  0          & {\rm cosh}\,t & {\rm sinh}\,t \\ 
              0          &  0          & {\rm sinh}\,t & {\rm cosh}\,t 
               \end{array} 
             \right] , \\ 
& P_{2, 1} :=\left[ 
             \begin{array}{cccc} 
              \cos \theta & 0 & -\sin \theta & 0 \\ 
               0          & 1 &   0          & 0 \\ 
              \sin \theta & 0 &  \cos \theta & 0 \\ 
               0          & 0 &   0          & 1 
               \end{array} 
             \right] , \ 
  P_{2, 2} :=\left[ 
             \begin{array}{cccc} 
              1          &   0           &   0           &   0           \\ 
              0          & {\rm cosh}\,t &   0           & {\rm sinh}\,t \\ 
              0          &   0           &   1           &   0           \\ 
              0          & {\rm sinh}\,t &   0           & {\rm cosh}\,t 
               \end{array} 
             \right] , \\ 
& P_{3, 1} :=\left[ 
             \begin{array}{cccc} 
               1          &  0          &   0          & 0 \\ 
               0          & \cos \theta & -\sin \theta & 0 \\ 
               0          & \sin \theta &  \cos \theta & 0 \\ 
               0          &  0          &   0          & 1 
               \end{array} 
             \right] , \ 
  P_{3, 2} :=\left[ 
             \begin{array}{cccc} 
              {\rm cosh}\,t & 0          & 0             & {\rm sinh}\,t \\ 
                0           & 1          & 0             &   0           \\ 
                0           & 0          & 1             &   0           \\ 
              {\rm sinh}\,t & 0          & 0             & {\rm cosh}\,t 
               \end{array} 
             \right] 
\end{split} 
\end{equation*} 
for $\theta , t\in \mbox{\boldmath{$R$}}$. 
Then the connected component $SO_0 (3,1)$ of the unit element $I_4$ 
of $SO(3,1)$ is generated by $P_{k, l}$ ($k=1, 2, 3$, $l=1, 2$, 
$\theta , t\in \mbox{\boldmath{$R$}}$). 
By direct computations, 
we have $\Tilde{T}_{-I_4} (\Theta_k )=\Theta_k$ for $k=1, 2, 3$, and 
\begin{equation*} 
\begin{split} 
&  [\Tilde{T}_{P_{1, 1}} (\Theta_1 ) \ 
    \Tilde{T}_{P_{1, 1}} (\Theta_2 ) \ 
    \Tilde{T}_{P_{1, 1}} (\Theta_3 )] 
  =[\Theta_1 \ \Theta_2 \ \Theta_3 ]Q_{1, 1} , \\ 
&  [\Tilde{T}_{P_{1, 2}} (\Theta_1 ) \ 
    \Tilde{T}_{P_{1, 2}} (\Theta_2 ) \ 
    \Tilde{T}_{P_{1, 2}} (\Theta_3 )] 
  =[\Theta_1 \ \Theta_2 \ \Theta_3 ]Q_{1, 2} , \\ 
& Q_{1, 1} :=\left[ 
    \begin{array}{ccc} 
     1 &  0          &   0          \\ 
     0 & \cos \theta & -\sin \theta \\ 
     0 & \sin \theta &  \cos \theta 
      \end{array} 
    \right] , \quad 
  Q_{1, 2} := 
    \left[ 
    \begin{array}{ccc} 
     1 &                 0        &                0        \\ 
     0 &            {\rm cosh}\,t & \sqrt{-1} {\rm sinh}\,t \\ 
     0 & -\sqrt{-1} {\rm sinh}\,t &           {\rm cosh}\,t 
      \end{array} 
    \right] , 
\end{split} 
\end{equation*} 
\begin{equation*} 
\begin{split} 
&  [\Tilde{T}_{P_{2, 1}} (\Theta_1 ) \ 
    \Tilde{T}_{P_{2, 1}} (\Theta_2 ) \ 
    \Tilde{T}_{P_{2, 1}} (\Theta_3 )] 
  =[\Theta_1 \ \Theta_2 \ \Theta_3 ]Q_{2, 1} , \\ 
&  [\Tilde{T}_{P_{2, 2}} (\Theta_1 ) \ 
    \Tilde{T}_{P_{2, 2}} (\Theta_2 ) \ 
    \Tilde{T}_{P_{2, 2}} (\Theta_3 )] 
  =[\Theta_1 \ \Theta_2 \ \Theta_3 ]Q_{2, 2} , \\ 
& Q_{2, 1} := 
    \left[ 
    \begin{array}{ccc} 
     \cos \theta & 0 & \sin \theta \\ 
      0          & 1 &  0          \\ 
    -\sin \theta & 0 & \cos \theta 
      \end{array} 
    \right] , \quad 
  Q_{2, 2} := 
    \left[ 
    \begin{array}{ccc} 
               {\rm cosh}\,t & 0 & \sqrt{-1} {\rm sinh}\,t \\ 
                    0        & 1 &                0        \\ 
    -\sqrt{-1} {\rm sinh}\,t & 0 &           {\rm cosh}\,t 
      \end{array} 
    \right] , 
\end{split} 
\end{equation*} 
and 
\begin{equation*} 
\begin{split} 
&  [\Tilde{T}_{P_{3, 1}} (\Theta_1 ) \ 
    \Tilde{T}_{P_{3, 1}} (\Theta_2 ) \ 
    \Tilde{T}_{P_{3, 1}} (\Theta_3 )] 
  =[\Theta_1 \ \Theta_2 \ \Theta_3 ]Q_{3, 1} , \\ 
&  [\Tilde{T}_{P_{3, 2}} (\Theta_1 ) \ 
    \Tilde{T}_{P_{3, 2}} (\Theta_2 ) \ 
    \Tilde{T}_{P_{3, 2}} (\Theta_3 )] 
  =[\Theta_1 \ \Theta_2 \ \Theta_3 ]Q_{3, 2} , \\ 
& Q_{3, 1} := 
    \left[ 
    \begin{array}{ccc} 
     \cos \theta & -\sin \theta & 0 \\ 
     \sin \theta &  \cos \theta & 0 \\ 
      0          &   0          & 1 
      \end{array} 
    \right] , \quad 
  Q_{3, 2} := 
    \left[ 
    \begin{array}{ccc} 
                {\rm cosh}\,t & \sqrt{-1} {\rm sinh}\,t & 0 \\ 
     -\sqrt{-1} {\rm sinh}\,t &           {\rm cosh}\,t & 0 \\ 
                     0        &                0        & 1 
      \end{array} 
    \right] . 
\end{split} 
\end{equation*} 
In particular, 
noticing $\Tilde{T}_{PP'} =\Tilde{T}_P \circ \Tilde{T}_{P'}$ 
for $P$, $P'\in SO(3, 1)$, 
we see that $\bigwedge^2_{\pm}\!X$ do not depend on the choice 
of $e=(e_1 , e_2 , e_3 , e_4 )$ as above. 
A Lie group $SO(3, \mbox{\boldmath{$C$}})$ is 
connected, of dimension $6$, and generated by 
$Q_{k, l}$ ($k=1, 2, 3$, $l=1, 2$, $\theta , t\in \mbox{\boldmath{$R$}}$). 
By $\Tilde{T}_{PP'} =\Tilde{T}_P \circ \Tilde{T}_{P'}$, 
we have 
a homomorphism $\Phi :SO_0 (3,1)\longrightarrow 
                      SO(3, \mbox{\boldmath{$C$}})$ 
satisfying $\Phi (P_{k, l} )=Q_{k, l}$ 
for $k=1, 2, 3$, $l=1, 2$, $\theta , t\in \mbox{\boldmath{$R$}}$. 
Since the differential of $\Phi$ is bijective, 
$\Phi :SO_0 (3,1)\longrightarrow SO(3, \mbox{\boldmath{$C$}})$ is 
a covering. 
In addition, $\Phi$ is injective, because 
\begin{itemize} 
\item{by the Cartan decomposition, 
$SL(2, \mbox{\boldmath{$C$}})$ is 
diffeomorphic to $SU(2)\times \mbox{\boldmath{$R$}}^3$,} 
\item{by a double covering from $SL(2, \mbox{\boldmath{$C$}})$ 
onto $SO_0 (3,1)$, 
$SO_0 (3,1)$ is diffeomorphic to $SO(3)\times \mbox{\boldmath{$R$}}^3$,} 
\item{by the Cartan decomposition, 
$SO(3, \mbox{\boldmath{$C$}})$ is also 
diffeomorphic to $SO(3)\times \mbox{\boldmath{$R$}}^3$.} 
\end{itemize} 
Therefore $SO_0 (3,1)$ is isomorphic to $SO(3, \mbox{\boldmath{$C$}})$ 
by $\Phi$. 
\end{rem} 

Let $E$ be an oriented vector bundle over a manifold $M$ of rank $4$. 
Let $h$ be a Lorentzian metric of $E$. 
The complexification $\bigwedge^2\!E\otimes \mbox{\boldmath{$C$}}$ 
of the two-fold exterior power $\bigwedge^2\!E$ of $E$ is decomposed into 
two subbundles $\bigwedge^2_{\pm}\!E$ of complex rank $3$ 
corresponding to the eigenvalues $\pm \sqrt{-1}$ respectively 
of Hodge's $*$-operator with respect to 
the Lorentzian metric $h$ and the orientation of $E$. 
The \textit{complex twistor spaces\/} $\hat{E}_{\pm}$ 
associated with $E$ are fiber bundles over $M$ 
such that the fibers $\hat{E}_{\pm , a}$ on each point $a$ of $M$ are 
given by 
\begin{equation*} 
\begin{split} 
  \{ z_1           \Theta_1 
    +z_2           \Theta_2 
    +z_3           \Theta_3  \ & | \ 
     z^2_1 +z^2_2 +z^2_3 =1, z_1 , z_2 , z_3 \in \mbox{\boldmath{$C$}} \} , 
  \\ 
  \{ z_1 \overline{\Theta_1} 
    +z_2 \overline{\Theta_2} 
    +z_3 \overline{\Theta_3} \ & | \ 
     z^2_1 +z^2_2 +z^2_3 =1, z_1 , z_2 , z_3 \in \mbox{\boldmath{$C$}} \} 
\end{split} 
\end{equation*} 
respectively 
for an ordered pseudo-orthonormal 
local frame field $(e_1 , e_2 , e_3 , e_4 )$ of $E$ 
defined on a neighborhood of $a$. 
In particular, 
$\Theta_1$, $\Theta_2$, $\Theta_3$ are local sections of $\hat{E}_+$ 
and $\overline{\Theta_1}$, $\overline{\Theta_2}$, $\overline{\Theta_3}$ 
are local sections of $\hat{E}_-$. 

Let $\nabla$ be an $h$-connection of $E$ 
and $\hat{\nabla}$ the induced connection of $\bigwedge^2\!E$ by $\nabla$. 
Then $\hat{\nabla}$ naturally gives a connection 
of $\bigwedge^2\!E\otimes \mbox{\boldmath{$C$}}$. 
We will prove 

\begin{pro}\label{pro:connL} 
The connection $\hat{\nabla}$ 
of $\bigwedge^2\!E\otimes \mbox{\boldmath{$C$}}$ gives 
connections of $\bigwedge^2_{\pm}\!E$. 
\end{pro} 

\vspace{3mm} 

\par\noindent 
\textit{Proof} \ 
Let $e=(e_1 , e_2 , e_3 , e_4 )$ be 
an ordered pseudo-orthonormal local frame field of $E$ 
giving the orientation and 
satisfying \eqref{pob} for the Lorentzian metric $h$ of $E$. 
Let $\omega =[\omega^i_j ]$ be the connection form of $\nabla$ 
with respect to $e$: $\nabla e=e\omega$. 
Then $\omega$ satisfies 
\begin{equation} 
\omega^k_l =-\omega^l_k \ (k, l=1, 2, 3), \quad 
\omega^k_4 = \omega^4_k \ (k=1, 2, 3),    \quad 
\omega^4_4 =0. 
\label{omega} 
\end{equation}  
By \eqref{omega}, we have 
\begin{equation} 
\begin{split} 
& [\hat{\nabla} \Theta_1 \ 
   \hat{\nabla} \Theta_2 \ 
   \hat{\nabla} \Theta_3 ] 
 =[\Theta_1 \ \Theta_2 \ \Theta_3 ]\hat{\omega} , \\ 
&  \hat{\omega} 
:= \left[ 
   \begin{array}{ccc} 
      0                               & 
    -\omega^3_2 +\sqrt{-1} \omega^4_1 & 
     \omega^3_1 +\sqrt{-1} \omega^4_2 \\ 
     \omega^3_2 -\sqrt{-1} \omega^4_1 & 
      0                               & 
    -\omega^2_1 +\sqrt{-1} \omega^4_3 \\ 
    -\omega^3_1 -\sqrt{-1} \omega^4_2 & 
     \omega^2_1 -\sqrt{-1} \omega^4_3 & 
      0 
     \end{array}    
   \right] . 
\end{split} 
\label{hatomega} 
\end{equation} 
Hence we obtain Proposition~\ref{pro:connL}. 
\hfill 
$\square$ 

\vspace{3mm} 

\begin{rem} 
The matrix $\hat{\omega}$ given in \eqref{hatomega} is 
a $1$-form valued in the Lie algebra of $SO(3, \mbox{\boldmath{$C$}})$. 
\end{rem} 

\section{Space-like surfaces in Lorentzian 4-dimensional \\ 
space forms}\label{sect:Ls} 

\setcounter{equation}{0} 

\subsection{The equations of Gauss, Codazzi and Ricci} 

Let $N$ be a 4-dimensional Lorentzian space form 
with constant sectional curvature $L_0$. 
If $L_0 =0$, 
then we can suppose $N=E^4_1$; 
if $L_0 >0$, 
then we can suppose $N=\{ x\in E^5_1 \ | \ \langle x, x\rangle =1/L_0 \}$; 
if $L_0 <0$, 
then we can suppose $N=\{ x\in E^5_2 \ | \ \langle x, x\rangle =1/L_0 \}$, 
where $\langle \ , \ \rangle$ is the metric of $E^5_1$ or $E^5_2$. 
In particular, 
if $L_0 = 1$, then $N$ is the      de Sitter $4$-space $S^4_1$; 
if $L_0 =-1$, then $N$ is the anti-de Sitter $4$-space $H^4_1$. 
Let $h$ be the Lorentzian metric of $N$. 
Let $M$ be a Riemann surface 
and $F:M\longrightarrow N$ a space-like and conformal immersion 
of $M$ into $N$. 
Let $(u, v)$, $g$, $\lambda$ be as in Section~\ref{sect:Riem}. 
Let $\Tilde{\nabla}$ denote the Levi-Civita connection 
of $E^4_1$, $E^5_1$ or $E^5_2$ according to $L_0 =0$, $>0$ or $<0$. 
Let $T_1$, $T_2$ be as in Section~\ref{sect:Riem}. 
Let $N_1$, $N_2$ be normal vector fields of $F$ 
satisfying \eqref{N1N2nt}. 
Then we have \eqref{dt1t2} with 
\begin{equation} 
\begin{split} 
S & =\left[ \begin{array}{ccccc} 
         \lambda_u         & \lambda_v & -\alpha_1  &  \beta_1   & 1 \\ 
        -\lambda_v         & \lambda_u & -\alpha_2  &  \beta_2   & 0 \\ 
         \alpha_1          & \alpha_2  &  \lambda_u &  \mu_1     & 0 \\ 
         \beta_1           & \beta_2   &  \mu_1     &  \lambda_u & 0 \\ 
         -L_0 e^{2\lambda} &  0        &   0        &   0        & 0 
              \end{array} 
     \right] , \\ 
T & =\left[ \begin{array}{ccccc} 
         \lambda_v & -\lambda_u        & -\alpha_2  &  \beta_2   & 0 \\ 
         \lambda_u &  \lambda_v        & -\alpha_3  &  \beta_3   & 1 \\ 
         \alpha_2  &  \alpha_3         &  \lambda_v &  \mu_2     & 0 \\ 
         \beta_2   &  \beta_3          &  \mu_2     &  \lambda_v & 0 \\ 
          0        & -L_0 e^{2\lambda} &   0        &   0        & 0 
              \end{array} 
     \right] , 
\end{split} 
\label{STLs} 
\end{equation} 
and $\alpha_k$, $\beta_k$ ($k=1, 2, 3$) and $\mu_l$ ($l=1, 2$) are 
real-valued functions. 
From \eqref{dt1t2}, 
we obtain $S_v -T_u =ST-TS$ with \eqref{STLs}. 
This is equivalent to 
the system of the equations of Gauss, Codazzi and Ricci. 
The equation of Gauss is given by 
\begin{equation} 
  \lambda_{uu} +\lambda_{vv} +L_0 e^{2\lambda} 
=-\alpha_1 \alpha_3 +\beta_1 \beta_3 +\alpha^2_2 -\beta^2_2 . 
\label{gaussLs}
\end{equation}
The equations of Codazzi are given by 
\begin{equation} 
\begin{split}
   (\alpha_1 )_v -(\alpha_2 )_u 
& = \alpha_2 \lambda_u +\alpha_3 \lambda_v 
   +\beta_2 \mu_1      -\beta_1  \mu_2  , \\ 
   (\alpha_2 )_v -(\alpha_3 )_u 
& =-\alpha_1 \lambda_u -\alpha_2 \lambda_v 
   +\beta_3  \mu_1     -\beta_2  \mu_2  , \\ 
   (\beta_1 )_v -(\beta_2 )_u 
& = \beta_2  \lambda_u +\beta_3  \lambda_v 
   +\alpha_2 \mu_1     -\alpha_1 \mu_2  , \\ 
   (\beta_2 )_v -(\beta_3 )_u 
& =-\beta_1  \lambda_u -\beta_2  \lambda_v 
   +\alpha_3 \mu_1     -\alpha_2 \mu_2 .  
\end{split} 
\label{codazziLs} 
\end{equation} 
The equation of Ricci is given by 
\begin{equation} 
  (\mu_1 )_v -(\mu_2 )_u = \alpha_1 \beta_2 -\alpha_2 \beta_1 
                          +\alpha_2 \beta_3 -\alpha_3 \beta_2 . 
\label{ricciLs} 
\end{equation}

Suppose that $N$ is oriented and 
that $(T_1 , T_2 , N_1 , N_2 )$ gives the orientation. 
Let $e_1$, $e_2$, $e_3$, $e_4$ be as in \eqref{e1e2e3e4}. 
The complex bundle $\bigwedge^2\!F^*T\!N \otimes \mbox{\boldmath{$C$}}$ is 
decomposed into 
two subbundles $\bigwedge^2_{\pm}\!F^*T\!N$ of complex rank $3$, 
and $\Theta_1$, $\Theta_2$, $\Theta_3$ as in \eqref{Omega+X} form 
a local frame field of $\bigwedge^2_+\!F^*T\!N$, 
and $\overline{\Theta_1}$, $\overline{\Theta_2}$, $\overline{\Theta_3}$ 
as in \eqref{Omega-X} form a local frame field of $\bigwedge^2_-\!F^*T\!N$. 
The \textit{complex twistor lifts\/} of $F$ are locally given 
by $\Theta_1$, $\overline{\Theta_1}$, 
which are determined by the immersion $F$ 
and depend on the choices of neither the isothermal coordinates $(u, v)$ 
nor normal vector fields $N_1$, $N_2$. 

Let $\nabla$ be as in Section~\ref{sect:Riem}. 
Then $\nabla$ induces a connection $\hat{\nabla}$ 
of $\bigwedge^2\!F^*T\!N \otimes \mbox{\boldmath{$C$}}$ naturally. 
In addition, by Proposition~\ref{pro:connL}, 
$\hat{\nabla}$ gives connections of $\bigwedge^2_{\pm}\!F^*T\!N$. 
By \eqref{dt1t2} with \eqref{STLs}, we obtain 
\begin{equation} 
\begin{split} 
    \hat{\nabla}_{T_1}  
   (\Theta_1 \ \Theta_2 \ \Theta_3 ) 
& =(\Theta_1 \ \Theta_2 \ \Theta_3 ) 
    \left[ \begin{array}{ccc} 
            0 & - W   & \sqrt{-1} Y \\ 
            W &   0   & \psi        \\ 
 -\sqrt{-1} Y & -\psi &  0 
             \end{array} 
    \right] , \\ 
    \hat{\nabla}_{T_2} 
   (\Theta_1 \ \Theta_2 \ \Theta_3 ) 
& =(\Theta_1 \ \Theta_2 \ \Theta_3 ) 
    \left[ \begin{array}{ccc} 
            0 & \sqrt{-1} Z   &   X   \\ 
 -\sqrt{-1} Z &           0   & -\phi \\ 
 -          X &          \phi &   0 
             \end{array} 
    \right] , 
\end{split} 
\label{hatnablaLs} 
\end{equation} 
where 
\begin{equation} 
\begin{array}{lcl} 
W:=\alpha_2 -\sqrt{-1} \beta_1  , & \ & X:=\alpha_2 +\sqrt{-1} \beta_3 , \\ 
Y:=\beta_2  -\sqrt{-1} \alpha_1 , & \ & Z:=\beta_2  +\sqrt{-1} \alpha_3  
\end{array} 
\label{WXYZLs} 
\end{equation} 
and 
\begin{equation} 
\phi :=\lambda_u -\sqrt{-1} \mu_2 , \quad 
\psi :=\lambda_v +\sqrt{-1} \mu_1 . 
\label{phipsi*} 
\end{equation} 
We have 
\begin{equation} 
W +\overline{W} =X +\overline{X} , \quad 
Y +\overline{Y} =Z +\overline{Z} . 
\label{W=XY=ZLs} 
\end{equation} 

Let $\hat{R}$ be the curvature tensor of $\hat{\nabla}$. 
Since $N$ is a space form of constant sectional curvature $L_0$, 
we obtain 
\begin{equation} 
  \hat{R} (T_1 , T_2 ) 
 (\Theta_1 \ \Theta_2 \ \Theta_3 ) 
=(\Theta_1 \ \Theta_2 \ \Theta_3 ) 
  \left[ \begin{array}{ccc} 
          0 &  0                & 0                \\ 
          0 &  0                & L_0 e^{2\lambda} \\ 
          0 & -L_0 e^{2\lambda} & 0 
           \end{array} 
  \right] . 
\label{hatRLs} 
\end{equation} 
On the other hand, 
we can compute the left side of \eqref{hatRLs} 
by \eqref{hatnablaLs}, and then we obtain  

\begin{pro}\label{pro:GCRLs} 
The functions $W$, $X$, $Y$, $Z$ as in \eqref{WXYZLs} satisfy 
not only \eqref{W=XY=ZLs} but also 
\begin{equation} 
 WX-YZ=L_0 e^{2\lambda} +\phi_u +\psi_v 
\label{gaussricciLs} 
\end{equation} 
and 
\begin{equation} 
\begin{split}
   Y_v +\sqrt{-1} X_u 
& =-\sqrt{-1} W\phi -Z\psi , \\ 
   W_v +\sqrt{-1} Z_u 
& =-\sqrt{-1} Y\phi -X\psi . 
\end{split} 
\label{codazziWXYZLs} 
\end{equation} 
\end{pro} 

\begin{rem} 
By direct computations as in \cite{AH}, 
we see that 
\eqref{gaussricciLs} is equivalent to 
the system of equations \eqref{gaussLs}, \eqref{ricciLs} and that 
the system of two equations in \eqref{codazziWXYZLs} is equivalent to 
the system of four equations in \eqref{codazziLs}. 
\end{rem} 

\subsection{Degeneracy of the complex twistor lifts} 

Let $K$, $\nabla^{\perp}$, $R^{\perp}$ be as in 
Section~\ref{sect:Riem}. 
The complex twistor lifts $\Theta_1$, $\overline{\Theta_1}$ of $F$ are 
said to be \textit{nondegenerate\/} (respectively, \textit{degenerate\/}) 
if  $\hat{\nabla}_{T_1} \Theta_1$ 
and $\hat{\nabla}_{T_2} \Theta_1$ are 
linearly independent (respectively, dependent) 
over $\mbox{\boldmath{$C$}}$ at each point of $M$. 
Then from \eqref{hatnablaLs} and \eqref{gaussricciLs}, 
we obtain an analogue of Proposition~\ref{pro:WX+YZ=0} 
for space-like surfaces in $N$. 

The complex twistor lifts $\Theta_1$, $\overline{\Theta_1}$ of $F$ are 
nondegenerate if and only if $\Delta :=WX-YZ \not= 0$. 
If we suppose $\Delta \not= 0$, 
then \eqref{codazziWXYZLs} can be rewritten into $(\phi , \psi )=(A, B)$, 
where 
\begin{equation} 
  \left[ \begin{array}{c} 
          A \\ 
          B 
           \end{array} 
  \right] 
:=\dfrac{1}{\Delta} 
  \left[ \begin{array}{cc} 
          \sqrt{-1} X & -\sqrt{-1} Z \\ 
                    Y & -          W 
           \end{array} 
  \right] 
  \left[ \begin{array}{c} 
          Y_v +\sqrt{-1} X_u \\ 
          W_v +\sqrt{-1} Z_u 
           \end{array} 
  \right] . 
\label{ABWXYZLs} 
\end{equation} 

Let $W$, $X$, $Y$, $Z$ be complex-valued functions 
of two variables $u$, $v$ satisfying \eqref{W=XY=ZLs}, 
$\Delta \not= 0$, 
$A_u +B_v =\Delta$ 
and $(A+\overline{A} )_v =(B+\overline{B} )_u$, 
where $A$, $B$ are functions constructed by $W$, $X$, $Y$, $Z$ 
as in \eqref{ABWXYZLs}. 
Then $W$, $X$, $Y$, $Z$, $\phi :=A$, $\psi :=B$ satisfy 
\eqref{gaussricciLs} with $L_0 =0$ and \eqref{codazziWXYZLs}. 
By $(A+\overline{A} )_v =(B+\overline{B} )_u$, 
there exists a real-valued function $\lambda$ satisfying 
\begin{equation} 
\lambda_u =\dfrac{1}{2} (A+\overline{A} ), \quad 
\lambda_v =\dfrac{1}{2} (B+\overline{B} ). 
\label{lambdaABLs} 
\end{equation} 
We set 
\begin{equation} 
\begin{array}{lclcl} 
\alpha_1 :=\dfrac{\sqrt{-1}}{2}  (Y-\overline{Y} ), & \ & 
\alpha_2 :=\dfrac{1}{2}          (X+\overline{X} ), & \ & 
\alpha_3 :=\dfrac{1}{2\sqrt{-1}} (Z-\overline{Z} ), \\ 
\ & \ & \ & \ & \ \\ 
\beta_1  :=\dfrac{\sqrt{-1}}{2}  (W-\overline{W} ), & \ & 
\beta_2  :=\dfrac{1}{2}          (Y+\overline{Y} ), & \ & 
\beta_3  :=\dfrac{1}{2\sqrt{-1}} (X-\overline{X} ) 
\end{array} 
\label{alphabetaLs} 
\end{equation} 
and 
\begin{equation} 
\mu_1 :=\dfrac{1}{2\sqrt{-1}} (B-\overline{B} ), \quad 
\mu_2 :=\dfrac{\sqrt{-1}}{2}  (A-\overline{A} ). 
\label{mu1mu2Ls} 
\end{equation} 
Then $\lambda$, $\alpha_k$, $\beta_k$ ($k=1, 2, 3$), $\mu_l$ ($l=1, 2$) 
satisfy \eqref{gaussLs} with $L_0 =0$, \eqref{codazziLs}, \eqref{ricciLs}. 
Therefore these functions determine an immersion $F$ into $E^4_1$ 
such that the complex twistor lifts of $F$ are nondegenerate, 
and this immersion is unique up to an isometry of $E^4_1$. 

Hence we obtain 

\begin{pro}\label{pro:nondegLsL=0} 
If the complex twistor lifts $\Theta_1$, $\overline{\Theta_1}$ 
of the immersion $F$ into $E^4_1$ are nondegenerate, 
then $\Delta \not= 0$, 
and $A$, $B$ as in \eqref{ABWXYZLs} satisfy 
$(A, B)=(\phi , \psi )$ with \eqref{phipsi*} 
and $A_u +B_v =\Delta$. 
In addition, 
complex valued functions $W$, $X$, $Y$, $Z$ 
of two variables $u$, $v$ satisfying \eqref{W=XY=ZLs}, 
$\Delta \not= 0$, $A_u +B_v =\Delta$ 
and $(A+\overline{A} )_v =(B+\overline{B} )_u$ give 
an immersion $F$ into $E^4_1$ 
such that the complex twistor lifts are nondegenerate, 
which is unique up to an isometry of $E^4_1$. 
\end{pro} 

Let $W$, $X$, $Y$, $Z$ be complex-valued functions 
of two variables satisfying \eqref{W=XY=ZLs}, 
$\Delta \not= 0$ and 
\begin{equation} 
\Delta -\overline{\Delta} =(A-\overline{A} )_u +(B-\overline{B} )_v , 
\label{DeltabarDelta} 
\end{equation} 
and suppose that 
\begin{equation} 
f:=\Delta -A_u -B_v \ 
 (=\overline{\Delta} -\overline{A}_u -\overline{B}_v ) 
\label{fLs} 
\end{equation} 
is nowhere zero, and satisfies 
\begin{equation} 
f_u =f(A+\overline{A} ), \quad 
f_v =f(B+\overline{B} ). 
\label{fufvLs} 
\end{equation} 
For a nonzero real number $L_0$, 
suppose $f/L_0 >0$ and set $\lambda :=(1/2)\log f/L_0$. 
Then $\lambda$ satisfies \eqref{lambdaABLs}. 
Let $\mu_1$, $\mu_2$ be functions as in \eqref{mu1mu2Ls}. 
Then $\phi$, $\psi$ as in \eqref{phipsi*} 
satisfy $\phi =A$, $\psi =B$. 
Therefore we obtain \eqref{gaussricciLs} and \eqref{codazziWXYZLs}, 
and for functions $\alpha_k$, $\beta_k$ ($k=1, 2, 3$) 
given as in \eqref{alphabetaLs}, 
we obtain \eqref{gaussLs}, \eqref{codazziLs}, \eqref{ricciLs}. 
Therefore the functions 
$\lambda$, $\alpha_k$, $\beta_k$ ($k=1, 2, 3$), $\mu_l$ ($l=1, 2$) 
determine an immersion $F$ into 
a 4-dimensional Lorentzian space form $N$ 
with constant sectional curvature $L_0 \not= 0$ 
such that the complex twistor lifts of $F$ are nondegenerate, 
and this immersion is unique up to an isometry of $N$. 

Hence we obtain 

\begin{pro}\label{pro:nondegLsLnot=0} 
If the complex twistor lifts $\Theta_1$, $\overline{\Theta_1}$ 
of the immersion $F$ into $N$ with $L_0 \not= 0$ are nondegenerate, 
then $\Delta \not= 0$, 
and $A$, $B$ as in \eqref{ABWXYZLs} satisfy 
$(A, B)=(\phi , \psi )$ with \eqref{phipsi*} 
and $f=L_0 e^{2\lambda}$ with \eqref{fLs}. 
In addition, 
a nonzero real number $L_0$ and 
complex-valued functions $W$, $X$, $Y$, $Z$ 
of two variables $u$, $v$ satisfying \eqref{W=XY=ZLs}, 
$\Delta \not= 0$, \eqref{DeltabarDelta} 
and \eqref{fufvLs} with \eqref{fLs} and $f/L_0 >0$ give 
an immersion $F$ into $N$ with constant sectional curvature $L_0$ 
such that the complex twistor lifts are nondegenerate, 
which is unique up to an isometry of $N$. 
\end{pro} 

\subsection{Horizontality of a complex twistor lift} 

Using \eqref{hatnablaLs}, 
we see that the covariant derivative of a complex twistor lift $\Theta_1$ 
by $\partial =\partial /\partial w$ is given by 
\begin{equation} 
 \hat{\nabla}_{\partial} \Theta_1 
=\dfrac{1}{2} ((W-Z)\Theta_2 +\sqrt{-1} (X-Y)\Theta_3 ). 
\label{dw} 
\end{equation} 
Whether $\hat{\nabla}_{\partial} \Theta_1$ vanishes or not 
is determined by the immersion $F$ and does not depend on 
the choice of a local complex coordinate $w$. 
By \eqref{dw}, $\hat{\nabla}_{\partial} \Theta_1$ vanishes 
if and only if 
\begin{equation} 
W=Z, \quad 
X=Y. 
\label{W=ZX=Y} 
\end{equation} 
Suppose that $F$ satisfies \eqref{W=ZX=Y}. 
Then applying \eqref{W=ZX=Y} to \eqref{gaussricciLs},  
we see that $F$ satisfies $K\equiv L_0$ and $R^{\perp} \equiv 0$, 
that is, 
the complex twistor lifts $\Theta_1$, $\overline{\Theta_1}$ are 
degenerate. 
Applying \eqref{W=ZX=Y} to \eqref{codazziWXYZLs}, we obtain 
\begin{equation} 
X_w =-W(\lambda +\gamma )_w , \quad 
W_w =-X(\lambda +\gamma )_w 
\label{XwWw} 
\end{equation} 
for a function $\gamma$ satisfying \eqref{gamma}. 
From \eqref{XwWw}, we obtain 
\begin{equation} 
((X+W)e^{\lambda  +\gamma} )_w =0, \quad 
((X-W)e^{-\lambda -\gamma} )_w =0. 
\label{XpmWw} 
\end{equation} 
By \eqref{XpmWw}, 
there exist 
anti-holomorphic functions        $\xi  =       \xi  (\overline{w} )$, 
                           $\Tilde{\xi} =\Tilde{\xi} (\overline{w} )$ 
of $w$ satisfying 
\begin{equation*} 
       \xi  (\overline{w} )=(X+W)e^{\lambda  +\gamma} , \quad 
\Tilde{\xi} (\overline{w} )=(X-W)e^{-\lambda -\gamma} . 
\end{equation*} 
Therefore we obtain 
\begin{equation*} 
\begin{split} 
W & =\dfrac{1}{2} 
   (        \xi  (\overline{w} )e^{-\lambda -\gamma} 
    -\Tilde{\xi} (\overline{w} )e^{\lambda  +\gamma} ), \\ 
X & =\dfrac{1}{2} 
   (        \xi  (\overline{w} )e^{-\lambda -\gamma}  
    +\Tilde{\xi} (\overline{w} )e^{\lambda  +\gamma} ). 
\end{split} 
\end{equation*} 
By \eqref{W=XY=ZLs}, there exists a real number $r$ 
such that $\Tilde{\xi}$ is constant so that $\Tilde{\xi} =\sqrt{-1} r$. 
Therefore we obtain 
\begin{equation} 
\begin{split} 
W & =\dfrac{1}{2} 
   ( \xi (\overline{w} )e^{-\lambda -\gamma} 
    -\sqrt{-1}      r   e^{\lambda  +\gamma} ), \\ 
X & =\dfrac{1}{2} 
   ( \xi (\overline{w} )e^{-\lambda -\gamma}  
    +\sqrt{-1}      r   e^{\lambda  +\gamma} ). 
\end{split} 
\label{WXw} 
\end{equation} 

Let $\lambda$ be a function of $u$, $v$ satisfying \eqref{K=L0}. 
Let $\gamma$ be a function of $u$, $v$. 
Let $\xi$ be an anti-holomorphic function of $w=u+\sqrt{-1} v$ 
and $r$ a real number. 
Let $W$, $X$ be as in \eqref{WXw}, 
and let $Y$, $Z$ be as in \eqref{W=ZX=Y}. 
Then $W$, $X$, $Y$, $Z$ and $\phi$, $\psi$ as in \eqref{phipsi*} 
with \eqref{gamma} satisfy 
\eqref{W=XY=ZLs}, \eqref{gaussricciLs} and \eqref{codazziWXYZLs}. 
Therefore, if we set 
\begin{equation} 
\begin{split} 
& \alpha_1 =-\beta_3 =\dfrac{1}{2\sqrt{-1}} (X-\overline{X} ), \\ 
& \alpha_2 = \beta_2 =\dfrac{1}{2}          (W+\overline{W} ), \\ 
& \alpha_3 =-\beta_1 =\dfrac{\sqrt{-1}}{2}  (W-\overline{W} ), 
\end{split} 
\label{alphabetaWX} 
\end{equation} 
then $\lambda$, $\alpha_k$, $\beta_k$ ($k=1, 2, 3$), $\mu_l$ ($l=1, 2$) 
satisfy \eqref{gaussLs}, \eqref{codazziLs}, \eqref{ricciLs}. 
Therefore these functions determine an immersion $F$ into $N$ 
such that the covariant derivative of a complex twistor lift $\Theta_1$ 
by $\partial =\partial /\partial w$ vanishes, 
and this immersion is unique up to an isometry of $N$. 

We can similarly study 
the condition $\hat{\nabla}_{\overline{\partial}} \Theta_1 =0$. 
This condition is equivalent to 
\begin{equation} 
W+Z=0, \quad 
X+Y=0. 
\label{W+Z=0X+Y=0} 
\end{equation} 
Suppose $\hat{\nabla}_{\overline{\partial}} \Theta_1 =0$. 
Then applying \eqref{W+Z=0X+Y=0} to \eqref{codazziWXYZLs}, we obtain 
\begin{equation*} 
X_{\overline{w}} =-W(\lambda -\gamma )_{\overline{w}} , \quad 
W_{\overline{w}} =-X(\lambda -\gamma )_{\overline{w}} , 
\label{XbarwWbarw} 
\end{equation*} 
instead of \eqref{XwWw}, and 
\begin{equation} 
\begin{split} 
W & =\dfrac{1}{2} 
    (\xi (w)e^{-\lambda +\gamma} +\sqrt{-1} re^{\lambda -\gamma} ), \\ 
X & =\dfrac{1}{2} 
    (\xi (w)e^{-\lambda +\gamma} -\sqrt{-1} re^{\lambda -\gamma} ) 
\end{split} 
\label{WX} 
\end{equation} 
for a holomorphic function $\xi =\xi (w)$ and a real number $r$, 
instead of \eqref{WXw}. 

Hence we obtain 

\begin{thm}\label{thm:deldelbar} 
The immersion $F$ satisfies 
either $\hat{\nabla}_{\partial} \Theta_1 =0$ 
or     $\hat{\nabla}_{\overline{\partial}} \Theta_1 =0$ 
if and only if $F$ satisfies \eqref{W=ZX=Y} or \eqref{W+Z=0X+Y=0}. 
If $F$ satisfies \eqref{W=ZX=Y} or \eqref{W+Z=0X+Y=0}, 
then 
\begin{itemize} 
\item[{\rm (a)}]{$K\equiv L_0$ and $R^{\perp} \equiv 0$,} 
\item[{\rm (b)}]{$W$, $X$ are represented as in 
either \eqref{WXw} or \eqref{WX}.} 
\end{itemize} 
In addition, 
for a function $\lambda$ satisfying \eqref{K=L0}, 
a function $\gamma$ and 
complex-valued functions $W$, $X$ represented as in \eqref{WXw} or \eqref{WX}, 
an immersion into $N$ satisfying 
either $\hat{\nabla}_{\partial} \Theta_1 =0$ 
or     $\hat{\nabla}_{\overline{\partial}} \Theta_1 =0$ is constructed 
by \eqref{W=ZX=Y} or \eqref{W+Z=0X+Y=0}, 
and unique up to an isometry of $N$. 
\end{thm} 

\begin{rem}\label{rem:deldelbar}  
By \eqref{WpmXpmcheck}, \eqref{WpmXpmcheck-}, 
   \eqref{WXw},         \eqref{WX},  
Theorem~\ref{thm:deldelbar} seems to be an analogue 
of Theorem~\ref{thm:pcsect}. 
\end{rem} 

\begin{rem} 
Suppose that the immersion $F$ 
satisfies $\hat{\nabla}_{\partial}            \Theta_1 =0$ 
or        $\hat{\nabla}_{\overline{\partial}} \Theta_1 =0$. 
Then the second fundamental form $\sigma$ satisfies 
the linearly dependent condition if and only if 
one of the following holds: 
\begin{itemize} 
\item[{\rm (i)}]{$\alpha_1 +\alpha_3 =0$, 
which means that $F$ has zero mean curvature vector;} 
\item[{\rm (ii)}]{$\alpha_1 =\alpha_3$ and $\alpha_2 =0$, 
which mean $\xi \equiv 0$.} 
\end{itemize} 
In any case, 
$\sigma$ satisfies the light-like linearly dependent condition, 
which is equivalent to the condition that 
a light-like normal vector field of $F$ is parallel (\cite{AH}). 
\end{rem} 

\begin{rem}\label{rem:cqd} 
Let $F:M\longrightarrow N$ be a space-like and conformal immersion. 
Let $w=u+\sqrt{-1} v$ be a local complex coordinate of $M$. 
Let $Q$ be a complex quartic differential on $M$ 
defined as in \eqref{QR}. 
Then we obtain 
\begin{equation*} 
  h(\sigma (\partial , \partial ), 
    \sigma (\partial , \partial )) 
  = \dfrac{e^{2\lambda}}{16} 
    (p_+ -2\sqrt{-1} q_+ )(p_- -2\sqrt{-1} q_- ) 
\end{equation*} 
for 
\begin{equation*} 
p_{\pm} :=\alpha_1 -\alpha_3 \pm (\beta_1 -\beta_3 ), \quad 
q_{\pm} :=\alpha_2           \pm  \beta_2 . 
\end{equation*} 
Therefore $Q$ is zero if and only if 
one of $p_{\pm} -2\sqrt{-1} q_{\pm}$ is zero. 
Noticing 
\begin{equation*} 
\begin{split} 
     p_+ -2\sqrt{-1} q_+ 
& =-\sqrt{-1} ( \overline{W} +\overline{X} +\overline{Y} +\overline{Z} ), \\ 
     p_- -2\sqrt{-1} q_- 
& =-\sqrt{-1} (           W  +          X  -          Y  -          Z  ), 
\end{split} 
\end{equation*} 
we see that if $F$ satisfies \eqref{W=ZX=Y} or \eqref{W+Z=0X+Y=0}, 
then $Q$ is zero. 
\end{rem} 

Let $F:M\longrightarrow N$ be a space-like and conformal immersion 
of $M$ into $N$ with zero mean curvature vector. 
Then we have $\alpha_3 =-\alpha_1$, $\beta_3 =-\beta_1$ 
and therefore $W=X$, $Y=Z$. 
Suppose that $\hat{\nabla}_{\partial} \Theta_1 =0$ vanishes. 
Then the second fundamental form $\sigma$ of $F$ satisfies 
the light-like linearly dependent condition. 
This is equivalent to a condition 
that the holomorphic quartic differential $Q$ 
defined as in \eqref{QR} vanishes (\cite{ando6}). 
Noticing $W=X$ and \eqref{WXw}, we obtain $r=0$, 
which means $e^{-\lambda -\gamma}  \xi (\overline{w} ) 
              =2(\alpha_2 -\sqrt{-1} \alpha_1 )$. 
Therefore, if $F$ has no totally geodesic points, 
then $\xi$ is nowhere zero, 
and then choosing a new local complex coordinate $w'=u'+\sqrt{-1} v'$ 
satisfying 
\begin{equation*} 
 \left( \dfrac{dw'}{dw} \right)^2 
=\dfrac{1-\delta \sqrt{-1}}{2} \overline{\xi (\overline{w})} 
\end{equation*} 
for $\delta \in \{ 1, -1\}$, 
we obtain 
\begin{equation*} 
\sigma (T'_1 , T'_2 )=\delta \sigma (T'_1 , T'_1 ), 
\end{equation*} 
where $T'_1 :=\partial /\partial u'$, $T'_2 :=\partial /\partial v'$. 
Since $\sigma$ is zero or light-like, 
we observe that $F$ is strictly isotropic in the sense of \cite{ando6}, 
that is, the mixed-type structure $\mu_F$ 
of the pull-back bundle $F^*T\!N$ given by $F$ 
satisfies $\mu_F \sigma (T'_1 , T'_1 )=\sigma (T'_1 , T'_2 )$ 
for $\delta =1$ or $-1$. 

Similarly, if $\hat{\nabla}_{\overline{\partial}} \Theta_1$ vanishes, 
then $r$ in \eqref{WX} is zero and in addition, 
if $F$ has no totally geodesic points, then $F$ is strictly isotropic. 

If we suppose 
that $\hat{\nabla}_{\partial} \Theta_1$ 
or   $\hat{\nabla}_{\overline{\partial}} \Theta_1$ vanishes 
so that $r$ in \eqref{WXw} or \eqref{WX} is zero, 
then we have $W=X$ and therefore $F$ has zero mean curvature vector. 

Hence we obtain 

\begin{cor}\label{cor:delbar} 
Suppose that 
the immersion $F$ satisfies $\hat{\nabla}_{\partial}            \Theta_1 =0$ 
or                          $\hat{\nabla}_{\overline{\partial}} \Theta_1 =0$. 
Then $F$ has zero mean curvature vector 
if and only if $r$ in either \eqref{WXw} or \eqref{WX} is zero. 
In addition, if $r=0$ and if $F$ has no totally geodesic points, 
then $F$ is strictly isotropic in the sense of\/ {\rm \cite{ando6}}. 
\end{cor}  

\begin{rem}\label{rem:ldcLs} 
Suppose that $F$ has zero mean curvature vector. 
Then $F$ satisfies 
either $\hat{\nabla}_{\partial}            \Theta_1 =0$ 
or     $\hat{\nabla}_{\overline{\partial}} \Theta_1 =0$ if and only if 
a light-like normal vector field of $F$ is parallel. 
Refer to \cite{ando3}, \cite{ando6}, \cite{AH} 
for characterizations of a space-like surface in $N$ 
with zero mean curvature vector 
and a parallel light-like normal vector field. 
The last assertion of Corollary~\ref{cor:delbar} is already shown in 
\cite{ando6}, and the above discussion directly gives it. 
\end{rem} 

\begin{rem} 
The conformal Gauss map of a space-like and Willmore surface 
in a $3$-dimensional Riemannian or Lorentzian space form 
with zero holomorphic quartic differential defines 
a space-like immersion into the      de Sitter $4$-space $S^4_1$ 
or                          the anti-de Sitter $4$-space $H^4_1$ 
with zero mean curvature vector, strict isotropicity 
and a parallel light-like normal vector field. 
See \cite{bryant2}, \cite{ando3} for such conformal Gauss maps. 
\end{rem} 

\section{Time-like surfaces in Lorentzian 4-dimensional \\ 
space forms}\label{sect:Lt} 

\setcounter{equation}{0} 

\subsection{The equations of Gauss, Codazzi and Ricci} 

Let $N$, $h$ be as in Section~\ref{sect:Ls}. 
Let $M$ be a Lorentz surface 
and $F:M\longrightarrow N$ a time-like and conformal immersion 
of $M$ into $N$. 
Let $(u, v)$, $g$, $\lambda$ be as in Section~\ref{sect:nt}. 
Let $\Tilde{\nabla}$, $T_1$, $T_2$ be as in Section~\ref{sect:Ls}. 
Let $N_1$, $N_2$ be normal vector fields of $F$ 
satisfying \eqref{N1N2}. 
Then we have \eqref{dt1t2} with 
\begin{equation} 
\begin{split} 
S & =\left[ \begin{array}{ccccc} 
       \lambda_u         & \lambda_v & -\alpha_1  & -\beta_1   & 1 \\ 
       \lambda_v         & \lambda_u &  \alpha_2  &  \beta_2   & 0 \\ 
       \alpha_1          & \alpha_2  &  \lambda_u & -\mu_1     & 0 \\ 
       \beta_1           & \beta_2   &  \mu_1     &  \lambda_u & 0 \\ 
       -L_0 e^{2\lambda} &  0        &   0        &   0        & 0 
              \end{array} 
     \right] , \\ 
T & =\left[ \begin{array}{ccccc} 
       \lambda_v &  \lambda_u        & -\alpha_2  & -\beta_2   & 0 \\ 
       \lambda_u &  \lambda_v        &  \alpha_3  &  \beta_3   & 1 \\ 
       \alpha_2  &  \alpha_3         &  \lambda_v & -\mu_2     & 0 \\ 
       \beta_2   &  \beta_3          &  \mu_2     &  \lambda_v & 0 \\ 
        0        &  L_0 e^{2\lambda} &   0        &   0        & 0 
              \end{array} 
    \right] , 
\end{split} 
\label{STLt} 
\end{equation} 
and $\alpha_k$, $\beta_k$ ($k=1, 2, 3$) and $\mu_l$ ($l=1, 2$) are 
real-valued functions. 
From \eqref{dt1t2}, 
we obtain $S_v -T_u =ST-TS$ with \eqref{STLt}. 
This is equivalent to 
the system of the equations of Gauss, Codazzi and Ricci. 
The equation of Gauss is given by 
\begin{equation} 
  \lambda_{uu} -\lambda_{vv} +L_0 e^{2\lambda} 
= \alpha_1 \alpha_3 +\beta_1 \beta_3 -\alpha^2_2 -\beta^2_2 . 
\label{gaussLt}
\end{equation}
The equations of Codazzi are given by 
\begin{equation} 
\begin{split}
   (\alpha_1 )_v -(\alpha_2 )_u 
& = \alpha_2 \lambda_u -\alpha_3 \lambda_v 
   -\beta_2 \mu_1      +\beta_1  \mu_2  , \\ 
   (\alpha_2 )_v -(\alpha_3 )_u 
& = \alpha_1 \lambda_u -\alpha_2 \lambda_v 
   -\beta_3  \mu_1     +\beta_2  \mu_2  , \\ 
   (\beta_1 )_v -(\beta_2 )_u 
& = \beta_2  \lambda_u -\beta_3  \lambda_v 
   +\alpha_2 \mu_1     -\alpha_1 \mu_2  , \\ 
   (\beta_2 )_v -(\beta_3 )_u 
& = \beta_1  \lambda_u -\beta_2  \lambda_v 
   +\alpha_3 \mu_1     -\alpha_2 \mu_2 .  
\end{split} 
\label{codazziLt} 
\end{equation} 
The equation of Ricci is given by 
\begin{equation} 
  (\mu_1 )_v -(\mu_2 )_u = \alpha_1 \beta_2 -\alpha_2 \beta_1 
                          -\alpha_2 \beta_3 +\alpha_3 \beta_2 . 
\label{ricciLt} 
\end{equation}

Suppose that $N$ is oriented and 
that $(T_1 , T_2 , N_1 , N_2 )$ gives the orientation. 
We set 
\begin{equation*} 
e_1 :=\dfrac{1}{e^{\lambda}} N_1 , \quad 
e_2 :=\dfrac{1}{e^{\lambda}} N_2 , \quad 
e_3 :=\dfrac{1}{e^{\lambda}} T_1 , \quad 
e_4 :=\dfrac{1}{e^{\lambda}} T_2 . 
\end{equation*} 
Then $\overline{\Theta_1}$, 
     $\overline{\Theta_2}$, 
     $\overline{\Theta_3}$ as in \eqref{Omega-X} 
form a local frame field of $\bigwedge^2_-\!F^*T\!N$. 
Let $\nabla$ be as in Section~\ref{sect:Riem}. 
Let $\hat{\nabla}$ be the connection 
of $\bigwedge^2\!F^*T\!N \otimes \mbox{\boldmath{$C$}}$ 
induced by $\nabla$. 
Then $\hat{\nabla}$ gives connections of $\bigwedge^2_{\pm}\!F^*T\!N$. 
By \eqref{dt1t2} with \eqref{STLt}, we obtain 
\begin{equation} 
\begin{split} 
    \hat{\nabla}_{T_1}  
   (\overline{\Theta_1} \ 
    \overline{\Theta_2} \ 
    \overline{\Theta_3} ) 
& =(\overline{\Theta_1} \ 
    \overline{\Theta_2} \ 
    \overline{\Theta_3} ) 
    \left[ \begin{array}{ccc} 
            0 & -\sqrt{-1}  W   & -\sqrt{-1}  Y   \\ 
  \sqrt{-1} W &             0   & -\sqrt{-1} \psi \\ 
  \sqrt{-1} Y &  \sqrt{-1} \psi &             0     
             \end{array} 
    \right] , \\ 
    \hat{\nabla}_{T_2} 
   (\overline{\Theta_1} \ 
    \overline{\Theta_2} \ 
    \overline{\Theta_3} ) 
& =(\overline{\Theta_1} \ 
    \overline{\Theta_2} \ 
    \overline{\Theta_3} ) 
    \left[ \begin{array}{ccc} 
            0 &            Z   & -           X   \\ 
           -Z &            0   & -\sqrt{-1} \phi \\  
            X & \sqrt{-1} \phi &             0 
             \end{array} 
    \right] , 
\end{split} 
\label{hatnablaLt} 
\end{equation} 
where 
\begin{equation} 
\begin{array}{lcl} 
W:=\alpha_2 +\sqrt{-1} \beta_1  , & \ & 
X:=\alpha_2 +\sqrt{-1} \beta_3 , \\ 
Y:=\beta_2  -\sqrt{-1} \alpha_1 , & \ & 
Z:=\beta_2  -\sqrt{-1} \alpha_3 
\end{array} 
\label{WXYZLt} 
\end{equation} 
and 
\begin{equation} 
\phi :=\lambda_u -\sqrt{-1} \mu_2 , \quad 
\psi :=\lambda_v -\sqrt{-1} \mu_1 . 
\label{phipsi*tL} 
\end{equation} 
We have \eqref{W=XY=ZLs}. 
Noticing that $N$ is a space form of constant sectional curvature $L_0$, 
we obtain 

\begin{pro}\label{pro:GCRLt} 
The functions $W$, $X$, $Y$, $Z$ satisfy 
not only \eqref{W=XY=ZLs} but also 
\begin{equation} 
WX +YZ +L_0 e^{2\lambda} +\phi_u -\psi_v =0 
\label{gaussricciLt} 
\end{equation} 
and 
\begin{equation} 
\begin{split}
   Y_v +\sqrt{-1} X_u 
& =-\sqrt{-1} W\phi -Z\psi , \\ 
   W_v -\sqrt{-1} Z_u 
& = \sqrt{-1} Y\phi -X\psi . 
\end{split} 
\label{codazziWXYZLt} 
\end{equation} 
\end{pro} 

\begin{rem} 
By direct computations as in \cite{AH}, 
we see that 
\eqref{gaussricciLt} is equivalent to 
the system of equations \eqref{gaussLt}, \eqref{ricciLt} and that 
the system of two equations in \eqref{codazziWXYZLt} is equivalent to 
the system of four equations in \eqref{codazziLt}. 
\end{rem} 

Let $K$, $\nabla^{\perp}$, $R^{\perp}$ be as in 
Section~\ref{sect:Riem}. 
Then from \eqref{hatnablaLt} and \eqref{gaussricciLt}, 
we obtain an analogue of Proposition~\ref{pro:WX+YZ=0} 
for time-like surfaces in $N$. 
In addition, noticing \eqref{gaussricciLt} and \eqref{codazziWXYZLt}, 
we obtain analogues of 
Propositions~\ref{pro:nondegLsL=0}, \ref{pro:nondegLsLnot=0}. 

\subsection{A complex twistor lift which is horizontal 
along integral curves of a light-like one-dimensional distribution} 

In the following, 
we study a time-like surface in $N$ 
such that the covariant derivative of $\overline{\Theta_1}$ 
along integral curves of a light-like one-dimensional distribution 
on the surface vanishes. 
By \eqref{hatnablaLt}, we have 
\begin{equation*} 
 \hat{\nabla}_{T_1 +T_2} \overline{\Theta_1} 
=(\sqrt{-1} W-Z)         \overline{\Theta_2} 
+(\sqrt{-1} Y+X)         \overline{\Theta_3} . 
\label{nabladu+dvT1} 
\end{equation*} 
Therefore $\hat{\nabla}_{T_1 +T_2} \overline{\Theta_1}$ vanishes 
if and only if $F$ satisfies 
\begin{equation} 
           (Z, Y) 
=\sqrt{-1} (W, X). 
\label{iWXZY} 
\end{equation} 
Applying \eqref{iWXZY} to \eqref{gaussricciLt}, 
we see that $F$ satisfies $K\equiv L_0$ and $R^{\perp} \equiv 0$. 
Applying \eqref{iWXZY} to \eqref{W=XY=ZLs}, 
we obtain $W=X$, and therefore 
\begin{equation} 
W=X=-\sqrt{-1} Y=-\sqrt{-1} Z. 
\label{WX-iY-iZ} 
\end{equation} 
Then by \eqref{codazziWXYZLt}, we obtain 
\begin{equation} 
   (\partial_u +\partial_v )X 
=-X(\partial_u +\partial_v )(\lambda -\sqrt{-1} \gamma ) , 
\label{dudvlgtL} 
\end{equation} 
which is equivalent to 
\begin{equation*} 
(\partial_u +\partial_v )(Xe^{\lambda -\sqrt{-1} \gamma} )=0. 
\label{dudvlgtL2} 
\end{equation*} 
Therefore there exists a complex-valued function $\xi$ 
of one variable satisfying 
\begin{equation} 
X=\xi (u-v)e^{-\lambda +\sqrt{-1} \gamma} . 
\label{W=X} 
\end{equation} 

Let $\lambda$ be a function of two variables $u$, $v$ 
satisfying \eqref{K=L0tn}. 
Let $\gamma$ be a function of $u$, $v$. 
Let $X$ be a function of $u$, $v$ 
represented as in \eqref{W=X} for a complex-valued function $\xi$ of $u-v$. 
Let $W$, $Y$, $Z$ be defined as in \eqref{WX-iY-iZ}. 
Let $\phi$, $\psi$ be defined as in \eqref{phipsi*tL} with \eqref{gamma}. 
Then these functions satisfy 
\eqref{W=XY=ZLs}, \eqref{gaussricciLt} and \eqref{codazziWXYZLt}. 
Therefore these functions give an immersion into $N$ 
satisfying $\hat{\nabla}_{T_1 +T_2} \overline{\Theta_1} =0$, 
which is unique up to an isometry of $N$. 

We can similarly study 
the condition $\hat{\nabla}_{T_1 -T_2} \overline{\Theta_1} =0$. 
This condition is equivalent to 
\begin{equation} 
W=X=\sqrt{-1} Y=\sqrt{-1} Z, 
\label{WXiYiZ} 
\end{equation} 
instead of \eqref{WX-iY-iZ}. 
Applying \eqref{WXiYiZ} to \eqref{codazziWXYZLt}, 
we obtain 
\begin{equation*} 
   (\partial_u -\partial_v )X 
=-X(\partial_u -\partial_v )(\lambda +\sqrt{-1} \gamma ) , 
\label{dudvlgtL-} 
\end{equation*} 
instead of \eqref{dudvlgtL}, and 
\begin{equation} 
X=\xi (u+v)e^{-\lambda -\sqrt{-1} \gamma} , 
\label{W=X-} 
\end{equation} 
instead of \eqref{W=X}. 

A time-like surface in $N$ satisfies 
\eqref{WX-iY-iZ} or \eqref{WXiYiZ} 
if and only if the surface has zero mean curvature vector 
so that the shape operator of any normal vector field is 
zero or light-like. 

Hence we obtain 

\begin{thm}\label{thm:ldcdztL} 
For a time-like and conformal immersion $F:M\longrightarrow N$, 
the following are mutually equivalent\/$:$ 
\begin{itemize} 
\item[{\rm (A)}]{the covariant derivatives of $\overline{\Theta_1}$ 
along integral curves of a light-like one-dimensional distribution 
on the surface vanish\/$;$} 
\item[{\rm (B)}]{$F$ has zero mean curvature vector so that 
the shape operator of any normal vector field is 
zero or light-like\/$;$} 
\item[{\rm (C)}]{$F$ satisfies \eqref{WX-iY-iZ} or \eqref{WXiYiZ}.} 
\end{itemize} 
If $F$ satisfies \eqref{WX-iY-iZ} or \eqref{WXiYiZ}, 
then 
\begin{itemize} 
\item[{\rm (a)}]{$K\equiv L_0$ and $R^{\perp} \equiv 0$,} 
\item[{\rm (b)}]{$X$ is represented as in 
either \eqref{W=X} or \eqref{W=X-}.} 
\end{itemize} 
In addition, 
for a function $\lambda$ satisfying \eqref{K=L0tn}, 
a function $\gamma$ of $u$, $v$, 
a function $\xi$ of one variable $u-v$ or $u+v$ and 
a function $X$ represented as in \eqref{W=X} or \eqref{W=X-}, 
an immersion into $N$ with the above {\rm (A)} is constructed 
by \eqref{WX-iY-iZ} or \eqref{WXiYiZ}, 
and unique up to an isometry of $N$. 
\end{thm} 

\begin{rem} 
Theorem~\ref{thm:ldcdztL} is an analogue of Theorem~\ref{thm:ldcdz}. 
\end{rem} 

\begin{rem} 
See \cite{AH} for characterizations of time-like surfaces 
with (B) in Theorem~\ref{thm:ldcdztL}, 
which can be obtained, based on \cite{ando6}. 
\end{rem} 

\begin{rem} 
Let $F:M\longrightarrow N$ be a time-like and conformal immersion. 
Let $\check{Q}$ be a paracomplex quartic differential on $M$ defined 
as in \eqref{QNt}. 
If $F$ satisfies (B) in Theorem~\ref{thm:ldcdztL}, 
then $\check{Q}$ is zero or null (see \cite{AH}). 
We can directly observe it, based on this situation. 
We have 
\begin{equation*} 
  h(\sigma (\check{\partial} , \check{\partial} ), 
    \sigma (\check{\partial} , \check{\partial} ))
  = \dfrac{e^{2\lambda}}{16} 
( ({\rm Re}\,\check{p} +2j{\rm Re}\,\check{q} )^2 
 +({\rm Im}\,\check{p} +2j{\rm Im}\,\check{q} )^2 ), 
\end{equation*} 
where $\check{p}$, $\check{q}$ are as in \eqref{checkpq}. 
Therefore $\check{Q}$ is zero or null if and only if 
\begin{equation} 
 ({\rm Re}\,(\check{p} +2\delta \check{q} ))^2 
+({\rm Im}\,(\check{p} +2\delta \check{q} ))^2 
=0 
\label{nullLt} 
\end{equation} 
for $\delta =1$ or $-1$. 
By \eqref{WXYZLt}, we have 
\begin{equation} 
\begin{split} 
  {\rm Re}\,(\check{p} +2\delta \check{q} ) 
& =\sqrt{-1} ( Y-\overline{Z} )+\delta (W+\overline{W} ), \\ 
  {\rm Im}\,(\check{p} +2\delta \check{q} ) 
& =\sqrt{-1} (-W+\overline{X} )+\delta (Y+\overline{Y} ) 
\end{split} 
\label{ReIm} 
\end{equation} 
for $\delta \in \{ 1, -1\}$. 
Therefore by \eqref{nullLt} and \eqref{ReIm}, 
if $F$ satisfies \eqref{WX-iY-iZ} or \eqref{WXiYiZ}, 
then $\check{Q}$ is zero or null. 
\end{rem}

\vspace{4mm} 

\par\noindent 
\footnotesize{Faculty of Advanced Science and Technology, 
              Kumamoto University \\ 
              2--39--1 Kurokami, Chuo-ku, Kumamoto 860--8555 Japan} 

\par\noindent  
\footnotesize{E-mail address: andonaoya@kumamoto-u.ac.jp} 

\end{document}